\documentclass[smallextended]{svjour3}

\usepackage{marvosym}

\usepackage[show]{ed}

\usepackage{tikz}
\usepackage{pgf}
\usetikzlibrary{arrows}
\usetikzlibrary{calc}
\usetikzlibrary{decorations.pathmorphing}

\usepackage{soul}
\setstcolor{red}
\setulcolor{red}

\usepackage{rotating} 

\usepackage{mathrsfs}
\DeclareMathAlphabet{\mathpzc}{OT1}{pzc}{m}{it}

\usepackage{color}
\usepackage{amsmath}
\usepackage{graphicx}
\usepackage{amsfonts}
\usepackage{amssymb}%
\usepackage{latexsym}
\usepackage{psfrag}
\usepackage{accents}

\spnewtheorem{theoremnew}{Theorem}{\bf}{}
\newtheorem{definitionnew}[theoremnew]{Definition}{\bf}{}
\spnewtheorem{lemmanew}[theoremnew]{Lemma}{\bf}{}
\spnewtheorem{assumptionnew}[theoremnew]{Assumption}{\bf}{}

\numberwithin{theoremnew}{section}
\numberwithin{equation}{section}
\numberwithin{figure}{section}

\spnewtheorem{remarknew}[theoremnew]{Remark}{\bf}{}
\newcommand{\noi}{\noindent}

\newcommand{\cC}{{\cal C}}

\newcommand{\cH}{{\cal H}}
\newcommand{\cI}{{\cal I}}
\newcommand{\cL}{{\cal L}}
\newcommand{\cP}{{\cal P}}

\newcommand{\cS}{{\cal S}}
\newcommand{\cT}{{\cal T}}

\newcommand{\bx}{\mathbf{x}}

\newcommand{\bxi}{\boldsymbol{\xi}}
\newcommand{\bn}{\mathbf{n}}

\newcommand{\ba}{\mathbf{a}}
\newcommand{\bbb}{\mathbf{b}}

\newcommand{\bv}{\mathbf{v}}

\newcommand{\by}{\mathbf{y}}

\newcommand{\tildu}{\widetilde{u}}

\newcommand{\re}{{\rm e}}
\newcommand{\ri}{{\rm i}}
\newcommand{\rd}{{\rm d}}

\newcommand{\HoDk}{H^1_k(\Omega_R)}
\newcommand{\HoDR}{H^1_R(\Omega_R)}

\newcommand{\beq}{\begin{equation}}
\newcommand{\eeq}{\end{equation}}
\newcommand{\beqs}{\begin{equation*}}
\newcommand{\eeqs}{\end{equation*}}
\newcommand{\bit}{\begin{itemize}}
\newcommand{\eit}{\end{itemize}}
\newcommand{\ben}{\begin{enumerate}}
\newcommand{\een}{\end{enumerate}}
\newcommand{\bal}{\begin{align}}
\newcommand{\eal}{\end{align}}
\newcommand{\bals}{\begin{align*}}
\newcommand{\eals}{\end{align*}}
\newcommand{\bse}{\begin{subequations}}
\newcommand{\ese}{\end{subequations}}
\newcommand{\bpr}{\begin{proposition}}
\newcommand{\epr}{\end{proposition}}
\newcommand{\bre}{\begin{remarknew}}
\newcommand{\ere}{\end{remarknew}}
\newcommand{\bpf}{\begin{proof}}
\newcommand{\epf}{\end{proof}}
\newcommand{\ble}{\begin{lemmanew}}
\newcommand{\ele}{\end{lemmanew}}
\newcommand{\bco}{\begin{corollary}}
\newcommand{\eco}{\end{corollary}}
\newcommand{\bex}{\begin{example}}
\newcommand{\eex}{\end{example}}
\newcommand{\bth}{\begin{theoremnew}}
\newcommand{\enth}{\end{theoremnew}}

\newcommand{\Rea}{\mathbb{R}}
\newcommand{\Com}{\mathbb{C}}

\newcommand{\Oi}{{\Omega_-}}

\newcommand{\Oe}{{\Omega_+}}

\newcommand{\GR}{{\Gamma_R}}

\newcommand{\eps}{\varepsilon}

\newcommand{\pdiff}[2]{\frac{\partial #1}{\partial #2}}

\newcommand{\gv}{\nabla v}

\newcommand{\tendi}{\rightarrow \infty}

\def\XXint#1#2#3{{\setbox0=\hbox{$#1{#2#3}{\int}$}
     \vcenter{\hbox{$#2#3$}}\kern-.5\wd0}}

\usepackage{hyperref}
\definecolor{myblue}{rgb}{0,0,0.6}
\hypersetup{colorlinks=true,
linkcolor=myblue,citecolor=myblue,filecolor=myblue,urlcolor=myblue}
\newcommand*{\N}[1]{\left\|#1\right\|}

\allowdisplaybreaks[4]

\newcommand{\tfa}{\text{ for all }}
\newcommand{\tfor}{\text{ for }}

\newcommand{\tin}{\text{ in }}
\newcommand{\ton}{\text{ on }}
\newcommand{\tas}{\text{ as }}
\newcommand{\tand}{\text{ and }}
\newcommand{\tst}{\text{ such that }}
\newcommand{\tfind}{\text{ find }}

\newcommand{\vertiii}[1]{{\left\vert\kern-0.25ex\left\vert\kern-0.25ex\left\vert #1
    \right\vert\kern-0.25ex\right\vert\kern-0.25ex\right\vert}}

\newcommand{\DtN}{{\rm DtN}_k}

\usepackage{soul}

\definecolor{jwcol}{RGB}{27, 137, 18}  

\definecolor{dalcol}{rgb}{0.8,0,0}

\definecolor{escol}{rgb}{0,0,0.8}
\definecolor{estcol}{rgb}{0,0.5,0}
\definecolor{esnewcol}{rgb}{0,0.5,0}

\newcommand{\hFEM}{h}

\newcommand{\supp}{{\rm supp}}

\newcommand{\Cosc}{C_{\rm{osc}}}

\newcommand{\scat}{{sc}}
\newcommand{\abs}[1]{{\lvert{#1}\rvert}}

\newcommand{\OR}{\Omega_R}

\newcommand{\truncbound}{\Gamma_R}
\newcommand{\Rscat}{R_{\rm sc}}
\newcommand{\Omegascat}{\Omega_{\rm sc}}

\newcommand{\OmegascataR}{\Omega_{{\rm sc}, R, \ba}}

\newcommand{\Omegascatb}{\Omega_{{\rm sc},\widetilde{\rho}, \bbb}}

\newcommand{\mymatrix}[1]{\mathsf{#1}}
\newcommand{\MA}{{\mymatrix{A}}}
\newcommand{\MI}{{\mymatrix{I}}}
\newcommand{\SPD}{{\mathsf{SPD}}}

\newcommand{\hatx}{\widehat{\bx}}

\newcommand{\CMS}{{C_{\rm MS}}}
\newcommand{\CPF}{{C_{\rm PF}}}
\newcommand{\Csol}{{C_{\rm sol}}}

\newcommand{\Creg}{{C_{H^2}}}
\newcommand{\Cmass}{{C_{\rm mass}}}
\newcommand{\Ccont}{{C_{\rm cont}}}
\newcommand{\Ccoer}{{C_{\rm coer}}}
\newcommand{\Cint}{{C_{\rm int}}}

\newcommand{\CTR}{{C_{\rm DtN}}}

\newcommand{\CR}{C}

\newcommand{\mythmname}[1]{(#1.)}

\begin{document}

\title{A sharp relative-error bound for the Helmholtz \lowercase{$h$}-FEM at high frequency}

\author{D.~Lafontaine \and E.~A.~Spence \and J.~Wunsch}

\date{\today}

\institute{D.~Lafontaine \at Department of Mathematical Sciences, University of Bath, Bath, BA2 7AY, UK, \email{D.Lafontaine@bath.ac.uk }
\and E.~A.~Spence \at Department of Mathematical Sciences, University of Bath, Bath, BA2 7AY, UK, \email{E.A.Spence@bath.ac.uk }
\and J.~Wunsch \at Department of Mathematics, Northwestern University, 2033 Sheridan Road, Evanston IL 60208-2730, US, \email{jwunsch@math.northwestern.edu}}

\maketitle

\begin{abstract}
For the $h$-finite-element method ($h$-FEM) applied to the Helmholtz equation, the question of how quickly the meshwidth $h$ must decrease with the frequency $k$ to maintain accuracy as $k$ increases has been studied since the mid 80's. Nevertheless, there still do not exist in the literature any $k$-explicit bounds on the \emph{relative error} of the FEM solution (the measure of the FEM error most often used in practical applications),
 apart from in one dimension.
The main result of this paper is the sharp result that, for the lowest fixed-order conforming FEM (with polynomial degree, $p$, equal to one), the condition ``$h^2 k^3$ sufficiently small" is sufficient for the relative error of the FEM solution in 2 or 3 dimensions to be controllably small (independent of $k$) for scattering of a plane wave by a nontrapping obstacle and/or a nontrapping inhomogeneous medium. 
We also prove relative-error bounds on the FEM solution  for arbitrary fixed-order methods applied to scattering by a nontrapping obstacle, but these bounds are not sharp for $p\geq 2$.
A key ingredient in our proofs is a result describing the oscillatory behaviour of the solution of the plane-wave scattering problem, which we prove using semiclassical defect measures.
\end{abstract}



\section{Introduction and informal statement of the main results}

\subsection{Introduction}\label{subsec:intro}
When solving the Helmholtz equation $\Delta u +k^2 u=0$ with the $h$ version of the finite-element method (where accuracy is increased by decreasing the meshwidth $h$ while keeping the polynomial degree $p$ constant), $h$ must decrease faster than $k^{-1}$ to maintain accuracy as $k$ increases; this is the so-called ``pollution effect" \cite{BaSa:00}.

A thorough investigation of how quickly $h$ must decrease with the frequency $k$ to maintain accuracy as $k$ increases was performed by Ihlenburg and Babu\v{s}ka in the mid 90's \cite{IhBa:95a,IhBa:97} on the 1-d model problem.
\beq\label{eq:1-d}
u'' +k^2 u = -f \quad \tin (0,1), \quad u(0)=0 \quad\tand \quad u'(1)-\ri k u(1)=0. 
\eeq
An explicit expression for the discrete Green's function for this problem is available, and Ihlenburg and Babu\v{s}ka used this to prove the following two sets of results:
\ben
\item The $h$-FEM is quasi-optimal in the $H^1$ semi-norm, with quasi-optimality constant independent of $k$, if $(hk^2/p)$ is sufficiently small; i.e.~there exists $c, C>0$, independent of $h, k,$ and $p$ such that, if $hk^2/p\leq c$, then 
\beqs
\N{\nabla (u- u_h)}_{L^2(0,1)} \leq C \min_{v_h\in \cH_h} \N{\nabla( u- v_h)}_{L^2(0,1)},
\eeqs
where $\cH_h$ is the appropriate conforming subspace of $H^1(0,1)$ of piecewise polynomials of degree $p$ on meshes of width $h$, and $u_h$ is the Galerkin solution; see \cite[Theorem 3]{IhBa:95a}, \cite[Theorem 4.13]{Ih:98}, \cite[Theorem 3.5]{IhBa:97} (when $p=1$ this result was proved earlier in \cite[Theorem 3.2]{AzKeSt:88}). 
The numerical experiments in \cite[Figures 8 and 9]{IhBa:95a} then indicated that, when $p=1$, the condition ``$hk^2$ sufficiently small" for quasi-optimality is necessary.
\item Under an assumption on the data $f$ (discussed below), 
the relative error in the $h$-FEM can be made arbitrarily small by, when $p=1$, making $hk^{3/2}$ sufficiently small and, when $p\geq 2$ and the data is sufficiently smooth (see \cite[Remark 4.28]{Ih:98}), making $h^{2p}k^{2p+1}$ sufficiently small.
 More precisely, \cite[Equation 3.25]{IhBa:95a}, \cite[Theorem 3.7]{IhBa:97}, \cite[Equation 4.5.15, \S4.6.4, and Theorem 4.27]{Ih:98} prove that there exists $C>0$, independent of $h$ and $k$ (but dependent on $p$) such that, if $hk$ is sufficiently small, then
the Galerkin solution $u_h$ exists and
\beq\label{eq:IhBa}
\frac{\N{u-u_h}_{H^1_k(0,1)}}{\N{u}_{H^1_k(0,1)}}\leq C \left(\left(\frac{hk}{p}\right)^p + k \left(\frac{hk}{p}\right)^{2p}\right),
\eeq
where the weighted $H^1$ norm $\|\cdot\|_{H^1_k(0,1)}$ is defined by \eqref{eq:1knorm} below. The numerical experiments in \cite[Figure 11]{IhBa:95a}, and \cite[Figure 4.13]{Ih:98} then indicated that, when $p=1$, the condition ``$h^{2}k^{3}$ sufficiently small" is necessary for the relative error to be bounded (in agreement with the earlier numerical experiments in \cite{BaGoTu:85} for small $k$).
\een

A note on terminology:~following \cite{IhBa:95a,IhBa:97,Ih:98}, we call the regime in $h,k,$ and $p$ where the solution is quasi-optimal (with constant independent of $k$) the 
\emph{asymptotic} regime, and the regime where the solution is not quasi-optimal the \emph{preasymptotic} regime.
For example, by the results in Points 1 and 2 above, when $p=1$ the asymptotic regime is when $hk^2$ is sufficiently small and 
the preasymptotic regime is when $hk^2 \gg 1$.

The (asymptotic) quasi-optimality results in Point 1 above have since been generalised to Helmholtz problems in 2 and 3 dimensions (and improved in the case $p \geq 2$). Indeed, the fact that the $h$-FEM with $p=1$ is quasi-optimal (with constant independent of $k$) in the full $H^1_k$ norm when $hk^2$ is sufficiently small was proved for the homogeneous Helmholtz equation on a bounded domain with impedance boundary conditions in \cite[Proposition 8.2.7]{Me:95} (in the case of constant coefficients) and \cite[Theorem 4.5 and Remark 4.6(ii)]{GrSa:20} (in the case of variable coefficients), and for scattering problems with variable coefficients in \cite[Theorem 3]{GaSpWu:20}. The fact that the $h$-FEM for $p \geq 2$ is quasi-optimal when $h^pk^{p+1}$ is sufficiently small was proved for a variety of constant coefficient Helmholtz problems  in \cite[Corollary 5.6]{MeSa:10}, \cite[Proof of Theorem 5.8]{MeSa:11}, and \cite[Theorem 5.1]{GaChNiTo:18}, and for a variety of problems including variable-coefficient Helmholtz problems in \cite[Theorem 2.15]{ChNi:20}; the condition ``$h^pk^{p+1}$ sufficiently small" is indicated to be sharp for quasi-optimality by, e.g., the numerical experiments in \cite[\S4.4]{ChNi:20}.

In contrast, the (preasymptotic) relative-error bound \eqref{eq:IhBa}  in Point 2 above has \emph{not} been obtained for any Helmholtz problem in 2 or 3 dimensions, even though numerical experiments indicate that the condition ``$h^{2p}k^{2p+1}$ sufficiently small" is necessary and sufficient for the relative error to be controllably small; see, e.g., \cite[Left-hand side of Figure 3]{DuWu:15}. The closest-available result is that, if $h^{2p}k^{2p+1}$ is sufficiently small, then
\beq\label{eq:Wu}
\N{u-u_h}_{H^1_k(D)}\leq C \left((hk)^p + k (hk)^{2p}\right)\N{f}_{L^2(D)},
\eeq
for the Helmholtz problem $\Delta u +k^2u =-f$ posed in a domain $D$ with either impedance boundary conditions on $\partial D$ or a perfectly matched layer (PML). Indeed, 
for the PML problem, \eqref{eq:Wu} is proved for $p=1$ in \cite[Theorem 4.4 and Remark 4.5(iv)]{LiWu:19} and \cite[Theorem 5.4]{GaChNiTo:18}. For the impedance problem, \eqref{eq:Wu} is proved for $p=1$ in \cite[Theorem 6.1]{Wu:14}, for $p \geq1$ in \cite[Corollary 5.2]{DuWu:15} (following earlier work by \cite{ZhWu:13}), and for $p\geq 1$ for the variable-coefficient Helmholtz equation $\nabla \cdot (\MA \nabla u) +k^2 nu =-f$ in \cite[\S2.3]{Pe:20} (under a nontrapping condition on $\MA$ and $n$). 

We highlight that, while \cite{DuWu:15}, \cite{GaChNiTo:18}, and \cite{LiWu:19} all prove results of the form \eqref{eq:Wu}, all the numerical experiments in these papers consider the \emph{relative error} (either in the $H^1$ norm 
\cite{DuWu:15}, \cite{LiWu:19}, or  the weighted $H^1$ norm \eqref{eq:1knorm} \cite{GaChNiTo:18}), illustrating that relative error is indeed the quantity of interest in practice.
An analogous situation is encountered in the preasymptotic error analyses of other Helmholtz FEMs in
\cite{FeXi:13,ZhDu:15,BuWuZh:16,DuZh:16,DuZh:17,WuZo:18,CaWu:20,DuWuZh:20,ZhWu:20}: all these papers prove bounds on the error in terms of the data, as in \eqref{eq:Wu}, but  
all the numerical experiments in these papers concerning the error consider the \emph{relative error}.








\subsection{The main results of this paper and their novelty}\label{sec:intro_main}
The two main results are the following:

\ben
\item[(a)] Theorem \ref{thm:main1} proves the relative-error bound \eqref{eq:IhBa} 
when $p=1$ for scattering 
of a plane wave by a nontrapping obstacle and/or a nontrapping inhomogeneous medium 
(modelled by
the PDE $\nabla \cdot (\MA \nabla u) +k^2 nu =0$ with variable $\MA$ and $n$)
 in 2 or 3 dimensions (see Definition \ref{def:planewave} below
 for the precise definition of the boundary-value problems considered). As highlighted above, the numerical experiments in 
\cite{BaGoTu:85,IhBa:95a,Ih:98} show that ``$h^2 k^3$ sufficiently small" is necessary for the relative error of the $h$-FEM with $p=1$ 
to be controllably small (independent of $k$), and so the result of Theorem \ref{thm:main1} is the sharp bound to which the title of the paper refers.
\item[(b)] Theorem \ref{thm:main2} proves for $p\geq 2$ a slightly-weaker bound than \eqref{eq:IhBa}, namely that 
\beq\label{eq:IhBa_alt}
\frac{\N{u-u_h}_{\HoDk}}{\N{u}_{\HoDk}}\leq C \big(hk + k (hk)^{p+1}
\big)
\eeq
for scattering of a plane wave by a nontrapping obstacle in 2 or 3 dimensions,
where $C$ in \eqref{eq:IhBa_alt} is independent of $h$ and $k$ but depends on $p$, with this dependence given explicitly in the theorem.
\een
As discussed above, these are the first-ever frequency-explicit relative-error bounds on the Helmholtz $h$-FEM in 2 or 3 dimensions. We recall the interest (highlighted at the end of the previous subsection) from \cite{FeXi:13,ZhWu:13,Wu:14,DuWu:15,ZhDu:15,BuWuZh:16,DuZh:16,DuZh:17,GaChNiTo:18,WuZo:18,LiWu:19,CaWu:20,DuWuZh:20,ZhWu:20} in proving such bounds. 

An additional novelty of Theorem \ref{thm:main1} is that it applies to the variable-coefficient Helmholtz equation, and all the constants in the relative-error bound are explicit, not only in $k$ and $h$, but also in the coefficients $\MA$ and $n$. The only other coefficient-explicit, preasymptotic FEM error bound on the variable-coefficient Helmholtz equation in the literature appears in  \cite[Theorem 2.39]{Pe:20}, where the bound \eqref{eq:Wu} is proved for the interior impedance problem when $h^{2p}k^{2p+1}$ is sufficiently small and $\MA$ and $n$ are nontrapping.
The only other coefficient-explicit FEM error bounds for the Helmholtz equation with variable $\MA$ and $n$ are in \cite{GrSa:20} and \cite{GaSpWu:20}. Both prove quasi-optimality under the condition ``$hk^2$ sufficient small" when $p=1$,
with \cite[Theorems 4.2 and 4.5]{GrSa:20} proving this result for the interior impedance problem and 
\cite[Theorem 3]{GaSpWu:20} proving this result for scattering by a nontrapping Dirichlet obstacle.

Our two main results, Theorems \ref{thm:main1} and \ref{thm:main2}, are proved 
for a particular class of Helmholtz problems, namely those corresponding to scattering by a plane wave, and not for the equation 
$\Delta u +k^2u =-f$ with general $f\in L^2$. We highlight that, for this latter class of problems, it is unreasonable to expect a relative-error bound such as \eqref{eq:IhBa} to hold, and thus the best one can do is prove bounds for a particular class of realistic data (as we do here). For example, consider the 1-d problem \eqref{eq:1-d} with 
\beq\label{eq:bad_data}
f(x) := -\big[ \exp(\ri k^n x) \chi(x)\big]'' - k^2\big[ \exp(\ri k^n x) \chi(x)\big],
\eeq
 where $\chi$ has compact support in $(0,1)$. The solution to \eqref{eq:1-d} is then $u(x) = \exp (\ri k^n x)\chi(x)$, which oscillates on a scale of $k^{-n}$, i.e., a smaller scale than $k^{-1}$ when $n>1$. The finite-element method with, say, $p=1$ and $hk^{3/2}$ small (and independent of $k$) will therefore not resolve this solution, and hence a bound such as \eqref{eq:IhBa} does not hold.  This example is nevertheless consistent with the previous results recalled in \S\ref{subsec:intro} since (i) the assumptions on the solution $u$ in \cite[First equation in \S3.4]{IhBa:95a} and \cite[Definition 3.2]{IhBa:97} exclude such data $f$, and (ii) 
with $f$ given by \eqref{eq:bad_data}, $\|f\|_{L^2(0,1)} \sim k^{2n}$ and $\|u\|_{H^1_k(0,1)}\sim k^n$, so that $\|f\|_{L^2(0,1)} \gg \|u\|_{H^1_k(0,1)}$, and 
the error estimate \eqref{eq:Wu} holds in this case because, although the absolute error on left-hand side of \eqref{eq:Wu} is large, the right-hand side of \eqref{eq:Wu} is larger.

\subsection{Discussion of these results in the context of using semiclassical analysis in the numerical analysis of the Helmholtz equation}\label{sec:featurenotdefect}

In the last $\sim$10 years, there has been growing interest in using results about the $k$-explicit analysis of the Helmholtz equation 
from \emph{semiclassical analysis} (a branch of \emph{microlocal analysis}) to design and analyse numerical methods for the Helmholtz equation\footnote{A closely-related activity is the design and analysis of numerical methods for the Helmholtz equation based on proving \emph{new} results about the $k\tendi$ asymptotics of Helmholtz solutions for polygonal obstacles; see \cite{ChLa:07,HeLaMe:13,HeLaCh:14,ChHeLaTw:15,He:15,GiChLaMo:19}.
}. The activity has so far occurred in, broadly speaking, five different directions:

\ben
\item The use of the results 
in \cite{MeTa:85} (on the rigorous $k\tendi$ asymptotics of the solution of the Helmholtz equation in the exterior of a smooth convex obstacle with strictly positive curvature) to design and analyse $k$-dependent approximation spaces for integral-equation formulations \cite{DoGrSm:07,GaHa:11,AsHu:14,EcOz:17,SoXiXi:17,LaBo:17,Ec:18,EcEr:19}, 
\item The use of the results in 
\cite{MeTa:85}, along with those in 
\cite{Ik:88} on scattering from several convex obstacles, to analyse algorithms for multiple scattering problems \cite{EcRe:09,AnBoEcRe:10,BoEcRe:17,EcAnBo:20}.
\item The use of bounds on the Helmholtz solution operator (also known as \emph{resolvent estimates}) due to 
\cite{Mo:75} and \cite{Va:75} (with the latter using the propagation of singularities results in 
\cite{MeSj:82}) to prove $k$-explicit bounds on both inverses of boundary-integral operators 
and the inf-sup constant of the domain-based variational formulation \cite{ChMo:08,Sp:14,BaSpWu:16,ChSpGiSm:20}, and also to analyse preconditioning strategies \cite{GaGrSp:15}.
\item The use of identities introduced in \cite{Mo:75} to prove coercivity of boundary-integral operators \cite{SpKaSm:15} and to introduce new coercive formulations of Helmholtz problems \cite{SpChGrSm:11,MoSp:14,GaMo:17a,DiMoSp:19,GaMo:19}.
\item The use of bounds on the restriction of quasimodes 
of the Laplacian to hypersurfaces from 
\cite{Tat,BGT,T,HTacy,christianson2014exterior,Ta:17} to prove sharp $k$-explicit bounds on boundary integral operators 
\cite{GaSm:15}, \cite[Appendix A]{HaTa:15}, \cite{Ga:19}, \cite{GaSp:18}, with these bounds then used to prove sharp $k$-explicit bounds on the number of  iterations when GMRES is applied to boundary-integral equations
\cite{GaMuSp:19}.
\een
The results of the present paper include a sixth direction. Namely, a key ingredient in our proofs of Theorems \ref{thm:main1} and \ref{thm:main2} (indeed, the ingredient that allows one to obtain a \emph{relative-error} bound instead of a bound in terms of the data, such as \eqref{eq:Wu}) is a result describing the oscillatory behaviour of the solution of the plane-wave scattering problem, which we prove using \emph{semiclassical defect measures}. 
These measures describe
where the mass in phase space of a Helmholtz solution is concentrated in the high-frequency limit (see the discussion in \S\ref{sec:defect1} below), and
were introduced in 
\cite{Ge:91} and \cite{LiPa:93}; see \cite{Bu:97} for more discussion on the history of defect measures.

\section{Formulation of the problem}\label{sec:form}

\begin{assumptionnew}[Assumptions on the domain and coefficients]\label{ass:1}

(i) $\Omega_- \subset\Rea^d, d=2,3,$ is a bounded open Lipschitz set such that its open complement $\Omega_+:= \Rea^d\setminus \overline{\Omega_-}$ is connected. 

(ii) $\MA \in C^{0,1}(\Omega_+ , \SPD)$
(where $\SPD$ is the set of $d\times d$ real, symmetric, positive-definite matrices)
is  such that $\supp(\MI- \MA)$ is compact in $\Rea^d$ and
there exist $0<A_{\min}\leq A_{\max}<\infty$ such that, 
for all $\bxi\in \Rea^d$,
\beq\label{eq:Alimits}
 A_{\min}|\bxi|^2 \leq 
\bxi^T\big(\MA(\bx)\bxi\big)
 \leq A_{\max}|\bxi|^2 \quad\text{ for almost every }\bx \in \Omega_+.
\eeq

(iii) $n\in L^\infty(\Omega_+,\Rea)$ is such that $\supp(1-n)$ is compact in $\Rea^d$ and there exist $0<n_{\min}\leq n_{\max}<\infty$ such that
\beq\label{eq:nlimits}
n_{\min} \leq n(\bx)\leq n_{\max}\quad \text{ for almost every } \bx \in \Omega_+.
\eeq
\end{assumptionnew}

\begin{figure}[h!]
\centering{
\includegraphics[width=0.6\textwidth]{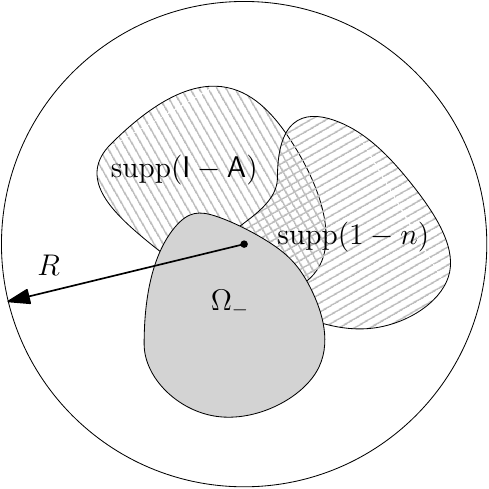}
}
\caption{A schematic of $\Omega_-$, the supports of $\MI-\MA$ and $1-n$, and $B_R$.}\label{fig:scatterer}
\end{figure}

Figure \ref{fig:scatterer} shows a schematic of $\Omega_-$ and the supports of $\MI-\MA$ and $1-n$.
Let the scatterer $\Omegascat$ be defined by $\Omegascat:= \Omega_-\cup \supp(\MI- \MA) \cup \supp(1-n)$ (i.e., the union of the shaded areas in Figure \ref{fig:scatterer}).
Given $R>0$ such that $\Omegascat \subset B_R$, where $B_R$ denotes the ball of radius $R$ about the origin, let $\Omega_R:= \Omega_+\cap B_R$.
Let $\Gamma_R:= \partial B_R$ and let $\Gamma:=\partial \Omega_-$. Let $\bn$ denote the outward-pointing unit normal vector field on both $\Gamma$ and $\Gamma_R$. We denote by $\partial_{\bn}$ the corresponding Neumann trace on $\Gamma$ or $\GR$ and $\partial_{\bn,\MA}$ the corresponding conormal-derivative trace.
We denote by $\gamma u$ the Dirichlet trace on $\Gamma$ or $\GR$.

\begin{definitionnew}[Helmholtz plane-wave scattering problem]\label{def:planewave}
Given $k>0$ and $\ba\in \Rea^d$ with $|\ba|=1$, let
$u^I(\bx) := \re^{\ri k \bx\cdot\ba}.$
Given $\Omega_-$, $\MA$, and $n$, as in Assumption \ref{ass:1},
we say $u\in H^1_{\rm{loc}}(\Omega_+)$ satisfies the \emph{Helmholtz plane-wave scattering problem} if 
\begin{align}
 \label{eq:PDE}
 \nabla\cdot(\MA\nabla u)+k^2 n u =0 \quad\tin\Omega_+,\quad
\text{either}\quad \gamma u =0 \quad\text{ or }\quad \partial_{\bn, \MA} u =0 \quad\ton \Gamma,
\end{align}
and $u^S:= u-u^I$ satisfies the Sommerfeld radiation condition 
\beq\label{eq:src}
\pdiff{u^S}{r}(\bx) - \ri k u^S(\bx) = o \left( \frac{1}{r^{(d-1)/2}}\right)
\eeq
as $r:= |\bx|\tendi$, uniformly in $\hatx:= \bx/r$. 
\end{definitionnew}

We call a solution of the Helmholtz equation satisfying the Sommerfeld radiation condition \eqref{eq:src} an \emph{outgoing} solution (so, in Definition \ref{def:planewave}, $u^S$ is outgoing).

 Define $\DtN: H^{1/2}(\truncbound) \rightarrow H^{-1/2}(\truncbound)$ to be the Dirichlet-to-Neumann map for the equation $\Delta u+k^2 u=0$ posed in the exterior of $B_R$ with the Sommerfeld radiation condition \eqref{eq:src}. 
When $\truncbound = \partial B_R$, for some $R>0$,  the definition of $\DtN$ in terms of Hankel functions and polar coordinates (when $d=2$)/spherical polar coordinates (when $d=3$) is given in, e.g., \cite[Equations 3.7 and 3.10]{MeSa:10}. 
Let 
\beqs
H_{0,D}^1(\OR):= \big\{ v\in H^1(\OR) : \gamma v=0 \ton \Gamma\big\}.
\eeqs
When Dirichlet boundary conditions are prescribed in \eqref{eq:PDE},  let 
\beq\label{eq:cHD}
\cH:= H_{0,D}^1(\OR);
\eeq
when Neumann boundary conditions are prescribed, let 
\beq\label{eq:cHN}
\cH:= H^1(\OR).
\eeq

\ble[Variational formulation of the Helmholtz plane-wave scattering problem]\label{lem:var}
With $u^I$, $\Omega_-$, $\MA$, $n$, $\OR$, and $\cH$ as above, 
define $\tildu\in \cH$ as the solution of the variational problem
\beq\label{eq:vf}
\text{ find } \tildu \in \cH \tst \quad a(\tildu,v)=F(v) \quad \tfa v\in \cH,
\eeq
where
\begin{align}\label{eq:sesqui}
a(\tildu,v)&:= \int_{\OR} 
\Big((\MA \nabla \tildu)\cdot\overline{\nabla v}
 - k^2 n \tildu\overline{v}\Big) - \big\langle \DtN (\gamma \tildu),\gamma v\big\rangle_{\GR},\quad\tand\quad\\
F(v)&:= \int_{\truncbound}\left(\partial_{\bn}u^I - \DtN(\gamma u^I)\right)\overline{\gamma v}.\nonumber
\end{align}
where $\langle\cdot,\cdot\rangle_{\GR}$ denotes the duality pairing on $\GR$ that is linear in the first argument and antilinear in the second.
Then $\tildu = u|_{\OR}$, where $u$ is the solution of the Helmholtz plane-wave scattering problem of Definition \ref{def:planewave}.
\ele

\noi For a proof of Lemma \ref{lem:var}, see, e.g., \cite[Lemma 3.3]{GrPeSp:19}. 
From here on we  denote the solution of the variational problem \eqref{eq:vf} by $u$, so that $u$ satisfies
\beq\label{eq:vf2}
\quad a(u,v)=F(v) \quad \tfa v\in \cH.
\eeq

\ble
\label{lem:wellposed}
The solution of the Helmholtz plane-wave scattering problem of Definition \ref{def:planewave} exists and is unique.
\ele

\bpf
Uniqueness follows from the unique continuation principle; see \cite[\S1]{GrPeSp:19}, \cite[\S2]{GrSa:20} and the references therein.
Since $a(\cdot,\cdot)$ satisfies a G\aa rding inequality (see \eqref{eq:Garding} below), 
Fredholm theory then gives existence. 
\epf

\paragraph{{\bf The $h$ finite-element method}}

Let $\cT_\hFEM$ be a family of triangulations of $\OR$ (in the sense of, e.g., \cite[Page 61]{Ci:91})  that is shape regular (see, e.g., \cite[Definition 4.4.13]{BrSc:08}, \cite[Page 128]{Ci:91}).
When Neumann boundary conditions are prescribed in \eqref{eq:PDE}, let
 \beq\label{eq:cHh}
 \cH_{\hFEM}:= \{v \in \CR(\overline{\Omega_R}) : v|_K \text{ is a polynomial of degree $p$ for each } K \in \cT_{\hFEM}\};
 \eeq
when Dirichlet boundary conditions are prescribed we impose the additional condition that elements of $\cH_{\hFEM}$ are zero on $\Gamma$; in both cases we then have $\cH_{\hFEM}\subset \cH$.
 The main results, Theorems \ref{thm:main1} and \ref{thm:main2} below
 require $\Gamma$ to be at least $C^{1,1}$. 
 For such $\OR$ it is not possible to fit $\partial \OR$ exactly with simplicial elements
(i.e.~when each element of $\cT_{\hFEM}$ is a simplex), and fitting $\partial \OR$ with isoparametric elements (see, e.g, \cite[Chapter VI]{Ci:91}) or curved elements (see, e.g., \cite{Be:89}) is impractical.
Some analysis of non-conforming error is therefore necessary, but
since this is very standard (see, e.g., \cite[Chapter 10]{BrSc:08}),
we ignore this issue here.

The second main result, Theorem \ref{thm:main2} (for $p \geq 2$ and analytic $\Gamma$), requires the triangulation $\cT_h$ to be quasi-uniform in the particular sense of 
\cite[Assumption 5.1]{MeSa:11}. Triangulations satisfying this assumption can be constructed by refining a fixed triangulation that has analytic element maps; see \cite[Remark 5.2]{MeSa:11}.

The finite-element method for the variational problem \eqref{eq:vf} is the Galerkin method applied to the variational problem \eqref{eq:vf}, i.e.
\beq\label{eq:FEM}
\text{ find } u_\hFEM \in \cH_\hFEM\tst\,\, a(u_\hFEM,v_\hFEM)=F(v_\hFEM) \,\, \tfa v_\hFEM\in  \cH_\hFEM.
\eeq
Observe that setting $v=v_h$ in \eqref{eq:vf2} and combining this with \eqref{eq:FEM} we obtain the \emph{Galerkin orthogonality} that
\beq\label{eq:GO}
a(u-u_h,v_h) =0 \quad\tfa v_h\in\cH_h.
\eeq

\section{Definitions of quantities involved in the statement of the main results}\label{sec:def}

Throughout the paper we assume that $R\geq R_0>0$ for some fixed $R_0>0$ and $k\geq k_0$ for some fixed $k_0>0$. For simplicity we assume throughout that
\beq\label{eq:assum}
k_0 R_0\geq 1 \quad\tand\quad hk\leq 1.
\eeq
Given a bounded open set $D$, we let the weighted $H^1$ norm, $\|\cdot\|_{H^1_k}$ be defined by
\begin{equation} \label{eq:1knorm}
\|u\|^2_{H^1_k(D)} :=\N{\nabla u}_{L^2(D)}^2 + k^2 \N{u}_{L^2(D)}^2.
\eeq

We now define quantities $\CTR_j, j=1,2, \Csol, \Cosc, \CPF,  \Creg, \Cint,$ and $\CMS$ that appear in the main results (Theorems \ref{thm:main1} and \ref{thm:main2}).
All of these are dimensionless quantities, independent of $k$, $h$, and $p$, but dependent on one or more of $\MA$, $n$, $\Omega_-$ (indicated below).

\paragraph{$\CTR_j$, $j=1,2$} By \cite[Lemma 3.3]{MeSa:10}, there exist $\CTR_j= \CTR_j(k_0 R_0)$, $j=1,2,$ such that 
\beq\label{eq:CDtN1}
\big|\big\langle \DtN(\gamma u), \gamma v\rangle_{\Gamma_R}\big\rangle\big| \leq \CTR_1 \N{u}_{H^1_k(\OR)}  \N{v}_{H^1_k(\OR)} 
\eeq
for all $u,v \in H^1(\OR)$ and for all $k\geq k_0$,
and 
\beq\label{eq:CDtN2}
- \Re \big\langle \DtN \phi,\phi\big\rangle_{\GR} \geq \CTR_2 R^{-1}\N{\phi}^2_{L^2(\GR)} 
\eeq
for all $\phi \in H^{1/2}(\GR)$ and for all $k\geq k_0.$

\paragraph{$\Csol$} 
We assume that $\MA$, $n$, and $\Omega_-$ are \emph{nontrapping} in the sense that there exists $\Csol =\Csol(\MA,n, \Oi, R,k_0)$ such that, given $f\in L^2(\Omega_R)$, 
the solution of the boundary value problem (BVP)
\begin{align*}
&\nabla\cdot(\MA\nabla v)+k^2 n v =-f \,\,\tin\Omega_+,\qquad
\text{either}\,\, \gamma v =0 \,\,\text{ or }\,\, \partial_{\bn, \MA} v =0 \,\,\ton \Gamma,
\end{align*}
and $v$ satisfies the Sommerfeld radiation condition \eqref{eq:src} (with $u^S$ replaced by $v$), satisfies the bound 
\beq\label{eq:Csol}
\N{v}_{\HoDk} \leq \Csol R\N{f}_{L^2(\Oe)} \quad \tfa k\geq k_0;
\eeq
observe that the factor $R$ on the right-hand side makes $\Csol$ dimensionless.
(Remark \ref{rem:trapping} discusses the situation where 
this nontrapping 
assumption is removed and $\Csol$ depends on $k$.)
This assumption holds if the obstacle $\Omega_-$ and the coefficients $\MA$ and $n$ are nontrapping in the sense that
all billiard trajectories (or, more precisely, Melrose--Sj\"ostrand
generalized bicharacteristics \cite[Section 24.3]{Ho:07}) starting in
an exterior neighbourhood of  $\Omega_-$ and evolving according to the
Hamiltonian flow defined by the symbol of \eqref{eq:PDE} escape from
that neighbourhood after some uniform time. For this flow to be
well-defined, $\Gamma$ must be $C^\infty$, and $\MA$ and $n$ must be
globally $C^{1,1}$ and $C^\infty$ in a neighbourhood of $\Gamma$;
note that the flow may in general be set-valued rather than unique in cases where the
  boundary is permitted to be infinite-order flat.
Assuming the uniqueness of the flow, an explicit expression for $\Csol$ in terms of $\MA, n,
\Oi,$ and $R$ is then given in  \cite[Theorems 1 and 2, and Equation
6.32]{GaSpWu:20}. However, the bound \eqref{eq:Csol} can be
established in situations with much less smoothness; indeed,
\cite[Theorems 2.5, 2.7, and 2.19]{GrPeSp:19} establishes
\eqref{eq:Csol} for a Dirichlet $C^0$ star-shaped obstacle and
$L^\infty$ $\MA$ and $n$ satisfying certain monotonicity assumptions. 
Furthermore, our arguments in the rest of the paper do not need the flow to be
well-defined on  
$\Omegascat:= \Omega_-\cup \supp(\MI- \MA) \cup \supp(1-n)$,
they only require that the bound \eqref{eq:Csol} holds. We can therefore define nontrapping in this weaker sense,
and work with scatterers of much lower smoothness than in standard microlocal-analysis settings.

\paragraph{$\Cosc$} 
By Theorem \ref{thm:Cosci} below, if $\MA$, $n$, and $\Omega_-$ are nontrapping then there exists
$\Cosc=\Cosc(\MA, n, \Omega_-)$ (`osc' standing for `oscillation') such that for $u$ a solution to the
Helmholtz plane-wave scattering problem of
Definition~\ref{def:planewave},
\beq\label{eq:Cosci}
\vert u\vert_{H^2(\Omega_R)}  \leq \Cosc k\N{u}_{H^1_k(\Omega_R)},
\eeq
where $\vert\cdot\vert_{H^2(\Omega_R)}$ denotes the $H^2$ semi-norm; i.e.~$\vert u\vert_{H^2(\Omega_R)}:= \sum_{|\alpha|=2}\int_{\Omega_R}|\partial^\alpha u|^2$.

\paragraph{$\CPF$}
By \cite[\S5.3]{BrSc:08}, \cite[Corollary A.15]{ToWi:05}, there exists $\CPF = \CPF(\Omega_-)$ (`PF' standing for `Poincar\'e--Friedrichs') such that
\beq\label{eq:CPF}
R^{-2} \N{v}^2_{L^2(\Omega_R)} \leq \CPF \left(R^{-1}\N{\gamma v}^2_{L^2(\GR)} + \N{\nabla v}^2_{L^2(\OR)}\right)
\eeq
for all $v\in H^1(\OR)$.

\paragraph{$\Creg$} By Theorem \ref{thm:Creg} below, there exists $\Creg=\Creg(\MA,  \Omega_-)$ such that, if $f\in L^2(\Omega_{R})$ and $v\in H^1(\Omega_{R})$ satisfy
\begin{subequations}\label{eq:H2equation}
\begin{align}
&\nabla \cdot (\MA \nabla v) = -f \,\, \tin \Omega_{R},\qquad\partial_\bn v = \DtN (\gamma v)\ton \GR, \tand\\
&\qquad\text{ either } \,\,\gamma v =0 \,\,\text{ or }\,\, \partial_\bn v =0 \quad\ton \Gamma,
\end{align}
\end{subequations}
then 
\beq\label{eq:Creg}
\big|v\big|_{H^2(\Omega_R)} \leq \Creg\left(
\N{f}_{L^2(\OR)} + R^{-1}\N{\nabla v}_{L^2(\OR)} +R^{-2} \N{v}_{L^2(\OR)}\right).
\eeq
The key point in \eqref{eq:Creg} is that, although $v$ in \eqref{eq:H2equation} depends on $k$ via the boundary condition on $\GR$, $\Creg$ is independent of $k$.

\paragraph{$\Cint$}
By, e.g., \cite[Equation 4.4.28]{BrSc:08}, \cite[Theorem 4.1]{ScZh:90} the \emph{nodal interpolant} $I_\hFEM : \CR(\overline{\Omega_R}) \rightarrow \cH_\hFEM$ is well-defined for functions in $H^2(\Omega_R)$ (for $d=2,3$) and satisfies
\beq\label{eq:Cint}
\N{v- I_\hFEM v}_{L^2(\Omega_R)} + \hFEM \N{\nabla(v- I_\hFEM v)}_{L^2(\Omega_R)} \leq \Cint \hFEM^2\vert v\vert_{H^2(\Omega_R)} ,
\eeq
for all $v\in H^2(\Omega_R)$, for some $\Cint$ that depends only on the shape-regularity constant of the mesh. 
As a consequence of \eqref{eq:Cint}, the definition of $\|\cdot\|_{\HoDk}$ \eqref{eq:1knorm}, and the assumption that $hk\leq 1$ \eqref{eq:assum}, we have
\beq\label{eq:Cint2}
\N{v-I_h v}_{\HoDk}\leq \sqrt{2} \Cint h |v|_{H^2(\OR)}.
\eeq

\paragraph{$\CMS$} By \cite[Lemma 3.4 and Proposition 5.3]{MeSa:11}
there exists $\CMS = \CMS(\Omega_-)$ (`MS' standing for
`Melenk--Sauter') so that, if $\Gamma$ is analytic, $\MA=\MI$, $n=1$,
and $\Omega_+$ is nontrapping, then the bound \eqref{eq:etabound2}
below holds. 

In \S\ref{sec:intro_main} we recalled that the only other frequency- and coefficient-explicit FEM error bounds for the variable-coefficient Helmholtz equation appear in  \cite[Theorems 4.2 and 4.5]{GrSa:20}, \cite[Theorem 3]{GaSpWu:20}, and \cite[Theorem 2.39]{Pe:20}. We note here that the constants in these bounds are expressed in terms of analogous quantities to those defined above.

\section{Statement and discussion of the main results}

\subsection{The main results}

The first theorem holds for any $p \geq1$, but is most relevant in the case $p=1$.
\bth\label{thm:main1} Let $u$ be the solution of the Helmholtz
    plane-wave scattering problem (Definition~\ref{def:planewave}).
Assume that both Assumption \ref{ass:1} and \eqref{eq:assum} hold, $\Omega_-$ is $C^{1,1}$, and $\MA,$ $n,$ and $\Omega_-$ are 
nontrapping. If
$p\geq 1$ and
\beq\label{eq:meshthreshold1}
h^2 k^3 \leq C_1,
\eeq
then the Galerkin solution $u_h$ to the variational problem \eqref{eq:FEM} exists, is unique, and satisfies the bound
\beq\label{eq:mainbound1}
\frac{\N{u-u_h}_{\HoDk}
}{
 \N{u}_{\HoDk}
}
 \leq C_2 hk + C_3 h^2 k^3,
\eeq
where 
\begin{align*}\nonumber
C_1 := &\frac{1}{4 (A_{\max} + \CTR_1)n_{\max}(\Creg)^2(\Cint)^2\Csol R
} 
\left( n_{\max} + \frac{1}{k_0 R_0 \Csol} + 2\right)^{-1}
\\ &\hspace{2cm}\times  
\bigg( 1 + \frac{\sqrt{2}}{
\min\big\{ \CTR_2 (\CPF)^{-1} \, ,\, A_{\min}(1+\CPF)^{-1}\big\}
}\bigg)^{-1},
\end{align*}
\beqs
C_2:= \frac{\sqrt{2} \Cint \Cosc}{A_{\min}}\big(\max\big\{A_{\max}, n_{\max}\big\} + \CTR_1\big),
\eeqs
and 
\begin{align*}
C_3&:= \frac{4\sqrt{2}}{\sqrt{A_{\min}}} \big( A_{\max} + \CTR_1\big)(\Cint)^2 \Creg\Csol R\Cosc \sqrt{n_{\max}+ A_{\min}}\\
&\hspace{4cm}\times
\left( n_{\max} + \frac{1}{k_0R_0 \Csol} + 2\right).
\end{align*}
\enth

\bth\label{thm:main2}
Let $u$ be the solution of the Helmholtz
    plane-wave scattering problem (Definition~\ref{def:planewave}).
Assume that both Assumption \ref{ass:1} and \eqref{eq:assum} hold, $\MA=\MI$, $n=1$, $\Omega_-$ is a nontrapping Dirichlet obstacle, $\Gamma$ is analytic, and the triangulation $\cT_h$ in the definition of $\cH_h$ \eqref{eq:cHh} satisfies the quasi-uniformity assumption \cite[Assumption 5.1]{MeSa:11}.
If
\beq\label{eq:meshthreshold2}
\frac{(hk)^2}{p} + \Csol R \frac{k(hk)^{p+1}}{p^p} \leq \widetilde{C}_1
\eeq
then the Galerkin solution $u_h$ to the variational problem \eqref{eq:FEM} exists, is unique, and satisfies the bound
\beq\label{eq:mainbound2}
\frac{
\N{u-u_h}_{\HoDk}
}{
 \N{u}_{\HoDk}
 }
 \leq \left(\widetilde{C}_2+ \frac{\widetilde{C}_3 \CMS}{p}\right) hk + \widetilde{C}_3 \CMS \Csol R \,\frac{k(hk)^{p+1}}{p^p},
\eeq
where
\beqs
\widetilde{C}_1 := \frac{1}{ 2\sqrt{2}(1+ \CTR_1) \Creg \CMS}\bigg( 1 + \frac{\sqrt{2}}{
\min\big\{ \CTR_2 (\CPF)^{-1}  ,(1+\CPF)^{-1}\big\}
}\bigg)^{-1},
\eeqs
\beqs
\widetilde{C}_2:= \sqrt{2} \Ccont \Cint \Cosc,\quad\tand\quad
\widetilde{C}_3:= 4\big(1 + \CTR_1\big) \Cint \Cosc .
\eeqs
\enth
Observe that (i) the condition \eqref{eq:meshthreshold2} is satisfied if $h^{p+1}k^{p+2}$ is sufficiently small, and (ii) the bound \eqref{eq:mainbound2} is of the form \eqref{eq:IhBa_alt}.

The result of Theorem \ref{thm:main2} might appear not to be a high-order result, since the lowest-order terms in \eqref{eq:meshthreshold2} and \eqref{eq:mainbound2} are $h^2$ and $h$, respectively. Nevertheless, for fixed $p$, if $k(hk)^{p+1}$ is sufficiently small, so that \eqref{eq:meshthreshold2} is satisfied, then 
\beqs
h \sim k^{-1-1/(p+1)} \quad\text{ so }\quad hk \sim k^{-1/(p+1)} \ll 1 \quad\tas k\tendi,
\eeqs
and the dominant term on the right-hand side of \eqref{eq:mainbound2} is that involving $k(hk)^{p+1}$.
We highlight that Theorem \ref{thm:main2}, along with the previous work discussed in \S\ref{subsec:intro}, shows that high-order methods suffer less from the pollution effect than low-order methods.

\subsection{How the main results are proved}\label{sec:idea1}

Theorems \ref{thm:main1} and \ref{thm:main2} are proved using the so-called \emph{elliptic-projection argument} or \emph{modified duality argument}, used to prove the bound \eqref{eq:Wu} on the solution in terms of the data. 
We first make some remarks about the history of this argument, and then outline our new contributions.

Recall that the classic duality argument, coming out of ideas introduced in \cite{Sc:74},
proves quasi-optimality of the Helmholtz FEM, and was used in, e.g., \cite{AzKeSt:88,IhBa:95a,Me:95,Sa:06,MeSa:10,MeSa:11,ChNi:18,ChNi:20,GaChNiTo:18,GrSa:20,GaSpWu:20}.
The elliptic-projection argument is a modification of this argument that allows one to prove results in the preasymptotic regime (as opposed to the asymptotic regime). The initial ideas were introduced in the Helmholtz context  in \cite{FeWu:09,FeWu:11} for interior-penalty discontinuous Galerkin methods, and 
then further developed for the standard FEM and continuous interior-penalty methods in \cite{Wu:14,ZhWu:13}. The argument has been subsequently used by \cite{DuWu:15,BaChGo:17,WuZo:18,ChNi:18,GaChNiTo:18,LiWu:19} (see, e.g., the literature review in \cite[\S2.3]{Pe:20}).

We note that \cite{FeWu:11} and \cite{Wu:14} also used an error-splitting argument (with this idea called ``stability-error iterative improvement" in these papers), and that error splitting ideas were also used in \cite{DuWu:15}, together with the idea of using discrete Sobolev norms in the duality argument. Although we do not use these ideas in this paper, one expects that they could be used to improve the $p$ dependence in Theorem \ref{thm:main2}, but see \cite[Remark 2.48]{Pe:20} for a discussion on the challenges in doing this.

Our three new contributions to the elliptic-projection argument are
(i) a rigorous proof, using semiclassical defect measures, of the bound \eqref{eq:Cosci} describing the oscillatory behaviour of the solution of the plane-wave scattering problem (see Theorem \ref{thm:Cosci} below), (ii) the proof of $H^2$ regularity, with constant independent of $k$, of the solution of Poisson's equation with the boundary condition $\partial_\bn v= \DtN(\gamma v)$ (see \eqref{eq:Creg} and Theorem \ref{thm:Creg}), and (iii) determining how all the constants in the elliptic-projection argument depend on $\MA, n, \Omega_-,$ and $R$. 

Regarding (i): oscillatory behaviour similar to \eqref{eq:Cosci} of Helmholtz solutions has been an assumption in many analyses of finite- and boundary-element methods; see, e.g., \cite[First equation in \S3.4]{IhBa:95a}, \cite[Definition 3.2]{IhBa:97}, \cite[Definition 4.6]{BuSa:06},
\cite[Definition 3.5]{BaSa:07}, \cite[Assumption 3.4]{DiMoSp:19}. However, to our knowledge, the only existing rigorous results proving such behaviour are \cite[Theorems 1.1 and 1.2]{GrLoMeSp:15} and \cite[Theorem 1.11(c)]{GaMuSp:19}. These results concern the Neumann trace of the solution of the Helmholtz plane-wave scattering problem with $\MA=\MI$ and $n=1$, and are then used in \cite{GrLoMeSp:15} and \cite{GaMuSp:19} to analyse boundary-element methods applied to this problem. In common with \eqref{eq:Cosci}, these results are obtained using semiclassical-analysis techniques.

Regarding (ii): the analogous result ($H^2$ regularity with constant independent of $k$) for Poisson's equation with the \emph{impedance boundary condition} $\partial_\bn v=\ri k \gamma v$ is central to the elliptic-projection argument for the Helmholtz equation with impedance boundary conditions. This result was explicitly assumed in \cite[Lemma 4.3]{FeWu:11}, implicitly assumed in \cite{Wu:14,ZhWu:13,BaChGo:17,ChNi:18},
 and recently proved in \cite{ChNiTo:20}. Our proof of \eqref{eq:Creg} uses (and makes $\MA$-explicit) arguments from \cite{ChNiTo:20}, which in turn use results from \cite{Gr:85}, adapting them to deal with the operator $\DtN$, instead of $\ri k$, in the boundary condition. 

Regarding (iii): while the standard duality argument applied to the Helmholtz equation discussed above has recently been made explicit in 
$\MA, n,$ and $\Omega_-$ in \cite{GrSa:20,GaSpWu:20} (as discussed in \S\ref{sec:intro_main}), 
the only places in the literature where the elliptic-projection argument is made explicit in $\MA, n,$ and $\Omega_-$ are 
the present paper and \cite[\S2.3]{Pe:20}, leading to the coefficient-explicit preasymptotic error bounds on the Helmholtz FEM at high-frequency in Theorem \ref{thm:main1} and \cite[Theorem 2.39]{Pe:20}. 
One area in which we expect these results to be applied is in the analysis of \emph{uncertainty quantification (UQ)} algorithms for the high-frequency Helmholtz equation with random coefficients, as discussed in the following remark.

\bre[The importance of coefficient-explicit FEM results for Helmholtz UQ]
To analyse UQ algorithms that use the standard Helmholtz variational formulation, one needs to understand how 
existence and uniqueness of the Galerkin solution is affected by the randomness in the coefficients.
One therefore needs coefficient-explicit existence and uniqueness results for the Galerkin solution for the Helmholtz equation with variable (deterministic) coefficients (such as in Theorem \ref{thm:main1} and \cite[Theorem 2.39]{Pe:20}); this issue is highlighted (but not fully analysed) in the analysis of Monte Carlo and Multi-level Monte Carlo methods in \cite[Chapter 5]{Pe:20}; see \cite[Assumption 5.1 and Remark 5.2]{Pe:20}.

The only other analyses of uncertainty quantification (UQ) algorithms for the high-frequency Helmholtz equation with random coefficients in the literature are  \cite{FeLiLo:15} and \cite{GaKuSl:20} (concerning Monte Carlo and Quasi-Monte Carlo methods, respectively). 
Because of the issue described in the previous paragraph, these papers use formulations of the Helmholtz equation where existence and uniqueness of the Galerkin solution is established 
for all $k, h, p,$ and for a class of (deterministic) coefficients
(\cite{FeLiLo:15} uses the 
interior-penalty discontinuous-Galerkin method of \cite{FeWu:09,FeWu:11} and \cite{GaKuSl:20} uses the coercive formulation of \cite{GaMo:19}).
This then ensures that the Galerkin solution exists and is unique for all realisations of the random coefficients; see the discussion at the beginning of \cite[\S4]{FeLiLo:15}. 
\ere

\subsection{Why does Theorem \ref{thm:main2} not cover scattering by an inhomogeneous medium?}
In both the elliptic-projection argument and the standard duality argument, 
a key role is played by the quantity $\eta(\cH_h)$ defined by \eqref{eq:eta} below, which describes how well solutions of the (adjoint of the) Helmholtz equation can be approximated in $\cH_h$.

In the case $p=1$ we estimate $\eta(\cH_h)$ using $H^2$ regularity of the solution (which holds when $\MA$ and $\Omega_-$ satisfy the assumptions of Theorem \ref{thm:main1}), leading to the bound \eqref{eq:etabound1} below.
When $p\geq 1$, $\MA=\MI$, $n=1$, $\Omega_-$ is a Dirichlet obstacle, and $\Gamma$ is analytic, \cite{MeSa:11} proved the bound \eqref{eq:etabound2} on $\eta(\cH_h)$, and we use this result to prove Theorem \ref{thm:main2}. The bound \eqref{eq:etabound2} was proved via a judicious splitting of the solution \cite[Theorem 4.20]{MeSa:11} into an analytic but oscillating part, and an $H^2$ part that behaves ``well" for large frequencies, and this splitting is only available for the exterior Dirichlet problem with $\MA=\MI$ and $n=1$.

We highlight that an alternative splitting procedure valid for Helmholtz problems with variable coefficients was recently developed in \cite{ChNi:20}, leading to an alternative proof of the bound on $\eta(\cH_h)$ \eqref{eq:etabound2} \cite[Lemma 2.13]{ChNi:20}.
However, this alternative procedure requires that $\DtN$ be approximated by $\ri k$ on $\GR$. Indeed, in \cite[Proof of Lemma 2.13]{ChNi:20} the solution is expanded in powers of $k$, i.e.~$u = \sum_{j=0}^{\infty} k^j u_j$, and then on $\GR$ one has $\partial_\bn u_{j+1} = \ri \gamma u_j$; this relationship between $u_{j+1}$ and $u_j$ on $\GR$ no longer holds if $\DtN$ is not approximated by $\ri k$.

\subsection{Approximating $\DtN$}
Implementing the operator $\DtN$ is computationally expensive, and so in practice one seeks to approximate this operator by \emph{either} imposing an absorbing boundary condition on $\GR$, \emph{or} using a PML. In this paper we follow the precedent established in \cite{MeSa:10,MeSa:11} of, when proving new results about the FEM for exterior Helmholtz problems, first assuming that $\DtN$ is realised exactly.
We remark, however, that if the two key ingredients in \S\ref{sec:idea1} (a proof of the oscillatory behaviour \eqref{eq:Cosci} and $H^2$-regularity, independent of $k$, of a Poisson problem) can be established when $\DtN$ is replaced by an absorbing boundary condition on $\GR$, then the result of Theorem \ref{thm:main1} carry over to this case. When an impedance boundary condition (i.e.~the simplest absorbing boundary condition) is imposed on $\GR$, the necessary Poisson $H^2$-regularity result is proved in \cite{ChNiTo:20}, but we discuss below in Remark \ref{rem:impedance} the difficulties in proving \eqref{eq:Cosci} in this case.

\subsection{Removing the nontrapping assumption}
\label{rem:trapping}
The only place in the proofs of Theorems \ref{thm:main1} and \ref{thm:main2} where the nontrapping assumption (i.e.~the fact that $\Csol$ in \eqref{eq:Csol} is independent of $k$) is used is in the proof of the bound \eqref{eq:Cosci} (in Theorem \ref{thm:Cosci} below).
We sketch in Remark \ref{rem:trapping2} below how \eqref{eq:Cosci} can be proved in the trapping case (i.e.~when $\Csol$ is not independent of $k$);
the rest of the proofs of Theorems \ref{thm:main1} and \ref{thm:main2} then go through as before. 
In the case of Theorem \ref{thm:main1}, the requirement for the relative error to be bounded independently of $k$ would then be that $h^2 k^3 \Csol$ be sufficiently small.
Under the strongest form of trapping, $\Csol$ can grow exponentially through a sequence of $k$s \cite[\S2.5]{BeChGrLaLi:11}, but is bounded polynomially in $k$ if a set of frequencies of arbitrarily-small measure is excluded \cite[Theorem 1.1]{LaSpWu:19}. However, it is not clear how sharp the requirement  ``$h^2 k^3 \Csol$ sufficiently small" for the relative error to be bounded is in these cases.

\section{Outline of the proof}
\label{sec:idea}
As highlighted in \S\ref{sec:idea1}, one of the novelties of this paper is that it makes the elliptic-projection argument explicit in the coefficients $\MA$ and $n$. 
However, this explicitness means that many of the expressions in the proofs are complicated (in the same way as the expressions in the results in Theorems \ref{thm:main1} and \ref{thm:main2} are complicated). In this section therefore, we give an outline of the proof, keeping track of the dependence on $k, h$, and $p$, but ignoring the dependence on $\MA, n, \Omega_-$, and $R$. We use the notation $a\lesssim b$ when $a\leq Cb$ with $C$ independent of $k, h,$ and $p$, but dependent on $\MA,n,\Omega_-,$ and $R$.

As in the standard duality argument coming out of ideas introduced in \cite{Sc:74} and then formalised in \cite{Sa:06}, our starting point is the fact that, since $a(\cdot,\cdot)$ satisfies the G\aa rding inequality \eqref{eq:Garding}, Galerkin orthogonality \eqref{eq:GO} and continuity of $a(\cdot,\cdot)$ \eqref{eq:continuity} imply that, for any $v_h\in \cH_h$,
\begin{align}\nonumber
&A_{\min} \N{u-u_h}^2_{\HoDk} \leq \Re a(u-u_h, u-v_h) + k^2 \big(n_{\max} + A_{\min}\big)\N{u-u_h}^2_{L^2(\OR)}\\
& \leq \Ccont \N{u-u_h}_{\HoDk} \N{u-v_h}_{\HoDk}
+ k^2  \big(n_{\max} + A_{\min}\big)\N{u-u_h}^2_{L^2(\OR)}.\label{eq:key1}
\end{align}

Recall (from, e.g., \cite[Theorem 2.5]{Sa:06}, \cite[Theorem 4.3]{MeSa:10}, \cite[Theorem 6.32]{Sp:15})
that the standard duality argument (related to the Aubin-Nitsche trick) shows that
\beq\label{eq:Schatz1}
\N{u-u_h}_{L^2(\OR)} \leq \Ccont \eta(\cH_h) \N{u-u_h}_{\HoDk},
\eeq
where $\eta(\cH_h)$, defined by \eqref{eq:eta} below, describes how well solutions of the adjoint problem are approximated in the space $\cH_h$.
Inputting \eqref{eq:Schatz1} into \eqref{eq:key1} one obtains quasi-optimality, with constant independent of $k$, if $k \eta(\cH_h)$ is sufficiently small.
Lemma \ref{lem:eta} below shows that $\eta(\cH_h)\lesssim h+(hk)^p$, and thus the condition ``$k \eta(\cH_h)$ sufficiently small" is satisfied if $h^pk^{p+1}$ is sufficiently small.

In contrast, the elliptic-projection argument, which we follow, shows that 
\beq\label{eq:AN2}
\N{u-u_h}_{L^2(\OR)} \lesssim \eta(\cH_h) \N{u-w_h}_{\HoDk} \quad \tfa w_h \in \cH_h,
\eeq
provided that $hk^2 \eta(\cH_h)$ is sufficiently small (see Lemma \ref{lem:AN} below).
Observe that \eqref{eq:AN2} is a stronger bound than \eqref{eq:Schatz1}, since $w_h$ on the right-hand side of \eqref{eq:AN2} is arbitrary.
The proof of \eqref{eq:AN2} in our setting of the plane-wave scattering problem requires the new Poisson $H^2$-regularity bound \eqref{eq:Creg}, which we prove in Theorem \ref{thm:Creg} below.

Inputting \eqref{eq:AN2} into \eqref{eq:key1}, choosing $w_h= v_h$, and using the inequality 
\beq\label{eq:Cauchy}
2\alpha\beta \leq \eps\alpha^2 + \eps^{-1}
\beta^2 \quad\tfa \alpha,\beta,\eps>0,
\eeq
on the first term on the right-hand side of \eqref{eq:key1}, we obtain that, if $hk^2 \eta(\cH_h)$ is sufficiently small, then,
for any $v_h\in \cH_h$,
\begin{align*}
\N{u-u_h}^2_{\HoDk} 
& \lesssim \big(1 + k^2 (\eta(\cH_h))^2 \big)\N{u-v_h}^2_{\HoDk};
\end{align*}
i.e.~quasi-optimality.
Assuming $H^2$ regularity of the solution, and using \eqref{eq:Cint2}, 
we obtain that, if $hk^2 \eta(\cH_h)$ is sufficiently small, then
\beq\label{eq:key2}
\N{u-u_h}_{\HoDk}^2 \lesssim \big(1  + k^2 (\eta(\cH_h))^2 \big)h^2 |u|^2_{H^2(\OR)}.
\eeq

In the standard elliptic-projection argument (see, e.g., \cite[\S5.5]{ChNi:18}) applied to the PDE $\Delta u + k^2 u=- f$,
an $H^2$-regularity bound similar to \eqref{eq:Csol} and the nontrapping bound \eqref{eq:Csol} are combined to give $|u|_{H^2(\OR)}\lesssim k \|f\|_{L^2(\OR)}$, and combining this with both \eqref{eq:key2} and the bound $\eta(\cH_h) \lesssim hk$ (see \eqref{eq:etabound1} below) proves the bound \eqref{eq:Wu} with $p=1$ on the Galerkin error in terms of the data when $h^2k^3$ is sufficiently small.

In contrast, in this paper we prove, using semiclassical defect measures,
 that the solution to the plane-wave scattering problem satisfies \eqref{eq:Cosci}, i.e.~$|u|_{H^2(\OR)}\lesssim k \N{u}_{\HoDk}$, (see Theorem \ref{thm:Cosci} below), and using this in \eqref{eq:key2}, along with the bounds on $\eta(\cH_h)$ in Lemma \ref{lem:eta}, we obtain the relative-error bounds \eqref{eq:mainbound1} 
and \eqref{eq:mainbound2}.

In summary, once one has proved the bound \eqref{eq:Cosci} (which we do via semiclassical analysis) 
and the Poisson $H^2$-regularity bound \eqref{eq:Creg} (which we do using results from \cite{Gr:85} and properties of $\DtN$),
if one ignores the technicalities of making the argument explicit in $\MA$, $n$, $\Omega_-$, and $R$, 
then the proof of a preasymptotic relative-error bound follows via a straightforward modification of the elliptic-projection argument.
Given the large and sustained interest (reviewed in \S\ref{subsec:intro}) in preasymptotic relative-error bounds for the Helmholtz FEM, we believe 
this fact  illustrates the advantage of approaching the numerical analysis of the Helmholtz equation from a perspective encompassing both numerical-analysis and semiclassical-analysis techniques. 

\section{Proof of the Poisson $H^2$-regularity result \eqref{eq:Creg}}

\bth\label{thm:Creg}
With $\MA$, $\Omega_-$, $\Gamma$, and $\Omega_R$ as in Section \ref{sec:form}, let $v\in H^1(\Omega_R)$ be the solution of the 
Poisson boundary value
problem \eqref{eq:H2equation}.
If $\Gamma$ is $C^{1,1}$, then $v\in H^2(\Omega_R)$ and the bound \eqref{eq:Creg} holds.
\enth

We follow the recent proof of the related regularity result  \cite[Theorem 3.1]{ChNiTo:20} (where $\DtN$ is replaced by $\ri k$, $\MA=\MI$, and $\Omega_-=\emptyset$) and start by recalling results from 
\cite{Gr:85}. 

\ble
\label{lem:grisvard1}Let $D$ be a bounded, convex, open set of $\mathbb{R}^{n}$
with $C^{2}$ boundary. Then, for all $\bv \in H^1(D; \Com^d)$,
\beq\label{eq:Grisvard}
\int_{D}\bigg(|\nabla\cdot\bv|^{2}-\sum_{i,j=1}^{n}\int_{D}\pdiff{v_i}{x_j}\overline{\pdiff{v_j}{x_i}}\,\bigg)
\geq-2\Re\big\langle(\gamma \bv)_{T},\nabla_{T}(\gamma \bv\cdot \bn)\big\rangle_{\partial D},
\eeq
where $\nabla_T$ is the surface gradient on $\partial D$ and $(\gamma \bv)_T:= \gamma \bv  - \bn(\gamma \bv \cdot \bn)$ is the tangential component of $\gamma \bv$.
\ele

\begin{proof}
The result with $\bv$ real follows from \cite[Theorem 3.1.1.1]{Gr:85} and the fact that 
the second fundamental form of $\partial D$ (defined in, e.g., \cite[\S3.1.1]{Gr:85}) is non-positive (see \cite[Proof of Theorem 3.1.2.3]{Gr:85}).
The result with $\bv$ complex follows in a straightforward way by repeating the argument in \cite[Theorem 3.1.1.1]{Gr:85} for complex $\bv$.
\end{proof}
\ble\mythmname{{\bf \cite[Lemma 3.1.3.4]{Gr:85}}}
\label{lem:grisvard2} If $\MA \in C^{0,1}(D,\SPD)$ satisfies \eqref{eq:Alimits} (with $\Omega_+$ replaced by $D$),
then, for all $v\in H^2(D)$,
\beq\label{eq:wacky0}
(A_{\min})^{2}\sum_{i,j=1}^{d}\left|\frac{\partial^{2}v}{\partial x_{i}\partial x_{j}}\right|^{2}\leq\sum_{i,j,\ell,m=1}^{d}A_{i \ell}A_{j m}\frac{\partial^{2}v}{\partial x_{j}\partial x_{\ell}}\frac{\partial^{2}\overline{v}}{\partial x_{i}\partial x_{m}}.
\eeq
\ele

As a first step to proving Theorem \ref{thm:Creg}, we prove it in the case when $\Omega_-=\emptyset$.

\ble\label{lem:Creg}
Let $\MA \in C^{0,1}(B_R,\SPD)$ satisfy \eqref{eq:Alimits} (with $\Omega_+$ replaced by $B_R$) and be such that $\supp(\MI-\MA) \subset\subset B_R$. 
Given $f\in L^{2}(B_{R})$, let $v\in H^{1}(B_{R})$ be the solution of 
\beq\label{eq:wacky1}
\nabla \cdot (\MA \nabla v) = -f \,\, \tin B_{R},\qquad\partial_\bn v = \DtN (\gamma v)\ton \GR.
\eeq
Then $v\in H^2(B_R)$ and 
\begin{align*}
&\vert v\vert_{H^{2}(B_R)} ^2 \leq \frac{2}{(A_{\min})^ 2} \bigg[\Vert f\Vert_{L^{2}(B_R)} ^2 + \bigg(d ^4 \Vert \nabla \MA\Vert_{L^{\infty}(B_R)} ^2\\
&\hspace{3.5cm}+ \frac{2}{(A_{\min})^2}d^8\Vert \MA\Vert_{L^{\infty}(B_R)}^{2} \Vert \nabla \MA\Vert_{L^{\infty}(B_R)}^{2} \bigg)\Vert \nabla v \Vert_{L^2(B_R)} ^ 2 \bigg],
\end{align*}
where $\nabla \MA$ denotes the derivative of $\MA$.
\ele

\begin{proof}
Let $w\in H^1(\Rea^d)$ be the outgoing solution of the following transmission problem 
\begin{align*}
\nabla\cdot (\MA \nabla w ) =- f \quad\tin B_R,\qquad \Delta w + k^2 w =0 \quad&\tin \Rea^d\setminus \overline{B_R},\\
\gamma w_+ = \gamma w_- \,\,\tand\,\, \partial_\bn w_+ = \partial_\bn w_- \quad&\ton \GR,
\end{align*}
where $w_-:= w|_{B_R}$ and $w_+ := w|_{\Rea^d\setminus B_R}$.
(Note that it is important here that $\MA=\MI$ in a neighbourhood of $\GR$, so that $\partial_{\bn,\MA} w_-= \partial_\bn w_-$.)
By the definition of the operator $\DtN$, $w_-= v$.
Since $\GR$ is $C^2$, the regularity result \cite[Theorem 5.2.1 and \S5.4b]{CoDaNi:10} implies that $w_- \in H^2(B_R)$
and $w_+ \in H^2_{\rm loc}(\Rea^d\setminus \overline{B_R})$;
 therefore $v\in H^2(B_R)$.

Since $v \in H^2(B_R)$ and $\MA$ is Lipschitz, $\MA\nabla v \in H^1(B_R)$ and we can apply Lemma \ref{lem:grisvard1} with $\bv:=\MA\nabla v$. Since $\MA=\MI$ near $\GR$, $\bv = \nabla v$ near $\GR$ and so the right-hand side of \eqref{eq:Grisvard} becomes
\beqs
-2\Re\big\langle \nabla_T (\gamma v), \nabla_T (\partial_{\bn} v) \big\rangle_{\Gamma_R}
=-2 \Re \big \langle \nabla_T (\gamma v), \nabla_T (\DtN (\gamma v))\big\rangle_{\Gamma_R},
\eeqs
where we have used the boundary condition in \eqref{eq:wacky1}. 

Now, $\DtN$ and $\nabla_T$ commute on $\GR$; this can be seen either by rotation invariance, or by using the definition of $\DtN$ and $\nabla_T$ in terms of Fourier series on $\GR$. Therefore, the inequality \eqref{eq:CDtN2}  implies that the right-hand side of \eqref{eq:Grisvard} is non-negative, hence
\beq\label{eq:wacky2}
\sum_{i,j,\ell,m=1}^{d}\int_{B_R}\frac{\partial}{\partial x_{j}}\left(A_{i\ell}\frac{\partial v}{\partial x_{\ell}}\right)\frac{\partial}{\partial x_{i}}\left(A_{j m}\frac{\partial \overline{v}}{\partial x_{m}}\right)\leq\Vert f\Vert_{L^{2}(B_R)}^2.
\eeq
The left-hand side of \eqref{eq:wacky2} equals
\begin{equation}
\sum_{i,j,\ell,m=1}^d\int_{\Omega}A_{i\ell}A_{j m}\frac{\partial^{2}v}{\partial x_{j}\partial x_{\ell}}\frac{\partial^{2}\overline{v}}{\partial x_{i}\partial x_{m}}+\sum_{i,j,\ell,m=1}^d\int_{\Omega}R_{i,j,\ell,m},\label{eq:wacky3}
\end{equation}
where
\begin{align*}
R_{i,j,\ell,m}&=\frac{\partial A_{i\ell}}{\partial x_{j}}\frac{\partial v}{\partial x_{\ell}}A_{j m}\frac{\partial^{2}\overline{v}}{\partial x_{i}\partial x_{m}}+A_{i\ell}\frac{\partial^{2}v}{\partial x_{j}\partial x_{\ell}}\frac{\partial A_{j m}}{\partial x_{i}}\frac{\partial \overline{v}}{\partial x_{m}}+\frac{\partial A_{i\ell}}{\partial x_{j}}\frac{\partial v}{\partial x_{\ell}}\frac{\partial A_{j m}}{\partial x_{i}}\frac{\partial \overline{v}}{\partial x_{m}}\\
&=:R_{i,j,\ell,m}^{1}+R_{i,j,\ell,m}^{2}+R_{i,j,\ell,m}^{3}.
\end{align*}
By the Cauchy-Schwarz inequality
\beqs
\left|\int_{B_R}R_{i,j,\ell,m}^{1}\right|+\left|\int_{B_R}R_{i,j,\ell,m}^{2}\right|\leq2\Vert \MA\Vert_{L^{\infty}(B_R)}\Vert \nabla  \MA\Vert_{L^{\infty}(B_R)}\Vert\nabla v\Vert_{L^{2}(B_R)}\vert v\vert_{H^2(B_R)}
\eeqs
and
\beqs
\left|\int_{B_R}R_{i,j,\ell,m}^{3}\right|\leq\Vert \nabla  \MA\Vert_{L^{\infty}(B_R)}^{2}\Vert\nabla v\Vert_{L^{2}(B_R)}^{2}.
\eeqs
We therefore obtain
\begin{align*}
\bigg|\sum_{i,j,\ell,m=1}^d\int_{B_R}R_{i,j,\ell,m}\bigg|
&\leq2d^{4}\Vert \MA\Vert_{L^{\infty}(B_R)}\Vert \nabla  \MA\Vert_{L^{\infty}(B_R)}\Vert\nabla v\Vert_{L^{2}(B_R)}\vert v\vert_{H^{2}(B_R)}\\
&\qquad+d^{4}\Vert \nabla  \MA\Vert_{L^{\infty}(B_R)}^{2}\Vert\nabla v\Vert_{L^{2}(B_R)}^{2}.
\end{align*}
Combining this with \eqref{eq:wacky0}, \eqref{eq:wacky2}, and \eqref{eq:wacky3}, we obtain 
\begin{align*}
&(A_{\min})^{2}\vert v\vert_{{H}^{2}(B_R)}^{2}\leq\Vert f\Vert_{L^{2}(B_R)}^{2}+2d^{4}\Vert \MA\Vert_{L^{\infty}(B_R)}\Vert \nabla  \MA\Vert_{L^{\infty}(B_R)}\Vert\nabla v\Vert_{L^{2}(B_R)}\vert v\vert_{H^{2}(B_R)}\\
&\hspace{4cm}+d^{4}\Vert \nabla  \MA\Vert_{L^{\infty}(B_R)}^{2}\Vert\nabla v\Vert_{L^{2}(B_R)}^{2}.
\end{align*}
Using \eqref{eq:Cauchy} on the second term on the right-hand side, we obtain the result.
\end{proof}

We now use Lemma \ref{lem:Creg} to prove Theorem \ref{thm:Creg}.

\bpf[Proof of Theorem \ref{thm:Creg}]
Let $0<R_0<R_1<R$ be such that $\overline{\Omega_-} \subset B_{R_0}$, and let $\chi\in C^{\infty}(\mathbb{R}^{d})$ be such that $0\leq\chi\leq1$ and
\begin{gather*}
\chi=0 \,\,\tin B_{R_0} \quad\tand \quad\chi = 1\,\, \tin \Rea^d\setminus\overline{B_{R_1}}.
\end{gather*}
We decompose $v$ as 
\beq\label{eq:u1u2}
v=\chi v+(1-\chi) v=:v_{1}+v_{2}.
\eeq
Then $v_{1} \in H^1(B_R)$ and satisfies
\beqs
\nabla\cdot(\MA\nabla v_{1})=
-\chi f +\nabla\chi\cdot(\MA \gv) + \nabla v \cdot(\MA\nabla \chi) + v\nabla\cdot(\MA \nabla\chi)
\quad\text{in }B_{R},
\eeqs
and $\partial_{\bn}v_{1}=\text{DtN}_{k}(\gamma v_{1})$ on $\GR$.
Lemma \ref{lem:Creg} implies that $v_1 \in H^2(B_R)$ and that there exists $C_4= C_4(\MA, d, \chi)>0$ such that
\beq\label{eq:reg_obst_1}
|v_1|_{H^2(\OR)}\leq C_4\left(\N{f}_{L^2(\OR)} + R^{-1}\N{\nabla v}_{L^2(\OR)} + R^{-2}\N{v}_{L^2(\OR)}\right),
\eeq
where (i) we have used the fact that $\nabla\chi =0$ in a neighbourhood of $\Omega_-$ to write all the norms as norms over $\OR$, and (ii) we have inserted the inverse powers of $R$ on the right-hand side to keep $C_4$ a dimensionless quantity.
On the other hand, $v_{2}$ satisfies 
\beqs
\nabla\cdot(\MA\nabla v_{2})=-(1-\chi)f -\nabla\chi\cdot(\MA \gv) - \nabla v \cdot(\MA\nabla \chi) - v\nabla\cdot(\MA \nabla\chi)
\quad\text{in }B_{R},
\eeqs
$v_2=0$ in $B_R\setminus B_{R_1}$, and \emph{either} $\gamma v_2=0$ \emph{or} $\partial_\bn v_2=0$ on $\Gamma$.

Since $\MA$ is Lipschitz, $A_{\min}>0$, and both $\Gamma$ and $\GR$ are $C^{1,1}$, 
\cite[Theorems 2.3.3.2, 2.4.2.5, and 2.4.2.7]{Gr:85} imply that, if $w\in H^1(\OR)$, $\nabla\cdot(\MA \nabla w)\in L^2(\OR)$, and \emph{either} $\gamma w= 0$ \emph{or} $\partial_\bn w=0$ on $\partial \OR$, then $w\in H^2(\Omega_-)$ and there exists $C_5=C_5(\MA, \Omega_-, d, R)>0$ such that
\begin{align*}
\vert w\vert_{H^{2}(\OR)}&\leq C_5\left(\N{\nabla\cdot(\MA\nabla w)-w}_{L^{2}(\OR)}+R^{-1}\N{\nabla w}_{L^2(\OR)}\right.\\
&\hspace{5.5cm}\left.+ R^{-2}\N{w}_{L^2(\OR)}
\right).
\end{align*}
Applying this with $w=v_2$, we obtain that 
\beq\label{eq:reg_obst_2}
|v_2|_{H^2(\OR)}\leq C_6\left(\N{f}_{L^2(\OR)} + R^{-1}\N{\nabla v}_{L^2(\OR)} +R^{-2}\N{v}_{L^2(\OR)}\right),
\eeq
and 
the bound 
\eqref{eq:Creg} follows from combining \eqref{eq:reg_obst_1} and \eqref{eq:reg_obst_2} using \eqref{eq:u1u2}.
\end{proof}

\section{The elliptic projection and associated results}\label{sec:ep}

Define the sesquilinear form $a_\star(\cdot,\cdot)$ by 
\beq\label{eq:astar}
a_\star(u,v) := \int_{\OR} \MA \nabla u \cdot\overline{\nabla v} - \big\langle \DtN \gamma u, \gamma v \big\rangle_{\GR}.
\eeq

Recall from \eqref{eq:cHD} and \eqref{eq:cHN} that $\cH$ equals \emph{either} $H_{0,D}^1(\OR)$ 
(with Dirichlet 
conditions in \eqref{eq:PDE}) \emph{or} $H^1(\OR)$ (with Neumann 
conditions).

\ble[Continuity and coercivity of $a_\star(\cdot,\cdot)$]\label{lem:astar1}
For all $u, v \in \cH$, 
\begin{align}\label{eq:astarcont}
\big|a_\star(u,v)\big|\leq \Ccont_\star \N{u}_{\HoDk} \N{v}_{\HoDk}
\,\,\tand
\,\,\Re a_\star(v,v) \geq \Ccoer_\star \N{v}^2_{\HoDR},
\end{align}
where
\beqs
\Ccont_\star:= A_{\max} + \CTR_1, \quad \Ccoer_\star:= 
\min\big\{ \CTR_2 (\CPF)^{-1} \, ,\, A_{\min}(1+\CPF)^{-1}\big\},
\eeqs
and
\beq\label{eq:HoDR}
\N{v}^2_{\HoDR}:= \N{\nabla u}^2_{L^2(\OR)} + \frac{1}{R^2} \N{v}^2_{L^2(\OR)}.
\eeq
\ele

\bpf
The first inequality in \eqref{eq:astarcont} follows from 
the inequality \eqref{eq:CDtN1} and the Cauchy--Schwarz inequality. The second inequality in \eqref{eq:astarcont} 
follows from \eqref{eq:CDtN2} and \eqref{eq:CPF}.
\epf

    As a consequence of Lemma \ref{lem:astar1}, we have 
\beq\label{eq:starnormequiv}
\Ccoer_\star \N{v}^2_{\HoDR} \leq \big|a_\star(v,v)\big| \leq \Ccont_\star \N{v}^2_{\HoDk} \quad\tfa v\in \cH,
\eeq
and we then define the new norm on $\cH$,
\beqs
\N{v}_\star := \sqrt{a_\star(v,v)}.
\eeqs

\ble[Bounds on the solution of the variational problem associated with $a_\star(\cdot,\cdot)$] 
The solution of the variational problem
\beqs
\tfind u \in \cH \tst a_\star(u,v) = (f,v)_{L^2(\OR)} \quad \tfa v\in \cH
\eeqs
satisfies 
\beq\label{eq:astarbound1}
\N{u}_{\HoDR} \leq \frac{R}{\Ccoer_\star} \N{f}_{L^2(\OR)}
\quad\tand\quad
\vert  u\vert_{H^2(\OR)} \leq \Creg_\star \N{f}_{L^2(\OR)},
\eeq
where
\beqs
\Creg_\star:=\Creg \left(1 + \sqrt{2}(\Ccoer_\star)^{-1}\right).
\eeqs
\ele

\bpf
Since $a_\star(\cdot,\cdot)$ is continuous and coercive in $\cH$, the first bound in \eqref{eq:astarbound1} follows from the Lax--Milgram theorem and the fact that
\beqs
\sup_{v\in \cH}\frac{\big| (f,v)_{L^2(\OR)}\big|}
{\N{v}_{\HoDR}} \leq R \N{f}_{L^2(\OR)},
\eeqs
by the definition of $\|\cdot\|_{\HoDR}$ \eqref{eq:HoDR}.
The second bound in \eqref{eq:astarbound1} follows from combining the first bound in \eqref{eq:astarbound1} and the bound \eqref{eq:Creg}.
\epf

We now define the particular Galerkin projection known in the literature as the ``elliptic projection" (see the discussion in \S\ref{sec:idea1}).

\begin{definitionnew}[Elliptic projection $\cP_h$]
Given $u\in \cH$, define $\cP_h u \in \cH_h$ by 
\beqs
a_\star(v_h, \cP_h u) = a_\star(v_h,u)\quad\tfa v_h\in \cH_h.
\eeqs
\end{definitionnew}
Since $a_\star(\cdot,\cdot)$ is continuous and coercive in $\cH$ by Lemma \ref{lem:astar1}, the Lax--Milgram theorem implies that $\cP_h$ is well defined. The definition of $\cP_h$ then immediately implies the Galerkin-orthogonality property that
\beq\label{eq:GO2}
a_\star(v_h ,u-\cP_h u) = 0 \quad\tfa v_h\in \cH_h.
\eeq

\ble[Approximation properties of $\cP_h$]
The elliptic projection $\cP_h$ satisfies
\begin{align}
\label{eq:ep1}
\N{u-\cP_h u}_\star &\leq \sqrt{\Ccont_\star} \min_{v_h\in \cH_h} \N{u-v_h}_{\HoDk} \qquad\tand\\
\label{eq:ep2}
\N{u-\cP_h u}_{L^2(\OR)} &\leq h \sqrt{2}\Cint \Creg_\star \sqrt{\Ccont_\star}  \N{u-\cP_h u}_\star
\end{align}
for all $u\in \cH$.
\ele
\bpf
By the Cauchy--Schwarz inequality $a_\star(\cdot,\cdot)$ is continuous
in the $\|\cdot\|_\star$ norm, and by definition,
$a_\star(\cdot,\cdot)$ is coercive in this norm. Therefore C\'ea's
lemma implies that
\beqs
\N{u-\cP_h u}_\star \leq \min_{v_h\in \cH_h} \N{u-v_h}_{\star},
\eeqs
and \eqref{eq:ep1} follows from the norm equivalence
\eqref{eq:starnormequiv}.

To prove \eqref{eq:ep2} we use the standard duality argument. Given $u\in \cH$, let $\xi$ be the solution of the variational problem
\beq\label{eq:astaradj}
\tfind \xi \in \cH \tst a_\star(\xi,v)= (u-\cP_h u,v)_{L^2(\OR)} \quad\tfa v\in \cH. 
\eeq
Then, by Galerkin orthogonality \eqref{eq:GO2} and continuity of $a_\star(\cdot,\cdot)$, for all $v_h\in \cH_h$,
\begin{align}\nonumber
\N{u-\cP_h u}^2_{L^2(\OR)} = a_\star( \xi ,u-\cP_h u ) &= a_\star( \xi-v_h ,u-\cP_h u )\\
&\leq \N{\xi-v_h}_\star \N{u-\cP_h u}_\star. \label{eq:astaradj2}
\end{align}
By the norm equivalence \eqref{eq:starnormequiv}, the consequence \eqref{eq:Cint2} of the definition of $\Cint$,
the definition of $\xi$ \eqref{eq:astaradj}, 
and the second bound in \eqref{eq:astarbound1},
\begin{align*}
\N{\xi-I_h \xi}_\star \leq\sqrt{\Ccont_\star} \N{\xi-I_h \xi}_{\HoDk} &\leq  \sqrt{\Ccont_\star}\sqrt{2}\Cint h |\xi|_{H^2(\OR)} ,\\
&\hspace{-1cm}\leq \sqrt{\Ccont_\star}  \sqrt{2} \Cint h \Creg_\star \N{u-\cP_h u}_{L^2(\OR)},
\end{align*}
and the result \eqref{eq:ep2} follows from combining this last inequality with \eqref{eq:astaradj2}.
\epf

\section{Adjoint approximability}

\begin{definitionnew}[Adjoint solution operator $\cS^*$]
Given $f\in L^2(\OR)$, let $\cS^*f$ be defined as the solution of the variational problem
\beq\label{eq:S*vp}
\tfind \cS^*f \in \cH \quad\tst\quad a(v, \cS^*f) = (v,f)_{L^2(\OR)} \quad\tfa v\in \cH.
\eeq
\end{definitionnew}
$\cS^*$ is therefore the solution operator of the adjoint problem to the variational problem \eqref{eq:vf} with data in $L^2(\OR)$.

Green's second identity applied to outgoing solutions of the Helmholtz equation implies that
$\big\langle \DtN \psi, \overline{\phi}\big\rangle_{\GR} =\big\langle \DtN \phi, \overline{\psi}\big\rangle_{\GR} $
(see, e.g., \cite[Lemma 6.13]{Sp:15}); thus $a(\overline{v},u) = a(\overline{u},v)$ and so the definition \eqref{eq:S*vp} implies that
\beq\label{eq:S*fkey}
a(\overline{\cS^*f}, v)= (\overline{f},v)_{L^2(\OR)}\quad\tfa v\in \cH;
\eeq
i.e.~$\cS^*f$ is the complex-conjugate of an outgoing Helmholtz solution.

Following \cite{Sa:06}, we define the quantity $\eta(\cH_h)$ by 
\beq\label{eq:eta}
\eta(\cH_h): = \sup_{f\in L^2(\OR)} \min_{v_h\in\cH_h} \frac{\N{\cS^*f- v_h}_{\HoDk}}{\N{f}_{L^2(\OR)}};
\eeq
observe that this definition implies that, given $f\in L^2(\OR)$, 
\beq\label{eq:eta2}
\text{ there exists } w_h \in \cH_H \tst \N{\cS^* f - w_h}_{\HoDk} \leq \eta(\cH_h) \N{f}_{L^2(\OR)}.
\eeq

\ble\label{lem:eta}
Assume that $\MA, n,$ and $\Omega_-$ are 
nontrapping (and so \eqref{eq:Csol} holds with $\Csol$ independent of $k$).

(i) If $\Gamma \in C^{1,1}$ and $\MA\in C^{0,1}$, then 
\beq\label{eq:etabound1}
\eta(\cH_h)\leq hk \left[ \sqrt{2} \Cint\Creg\Csol R\left( n_{\max} + \frac{1}{k_0 R_0 \Csol} + 2\right) \right].
\eeq
(ii) If $\Omega_-$ is a Dirichlet obstacle (so that $\cH= H^1_{0,D}(\Omega_R)$), $\Gamma$ is analytic, $\MA= \MI$, $n=1$, 
and the triangulation $\cT_h$ in the definition of $\cH_h$ \eqref{eq:cHh} satisfies the quasi-uniformity assumption \cite[Assumption 5.1]{MeSa:11}, then there exists $\CMS = \CMS(\Omega_-)$ such that
\beq\label{eq:etabound2}
\eta(\cH_h) \leq \CMS \left[ \frac{h}{p} + \Csol R\left(\frac{hk}{p}\right)^p\right].
\eeq
\ele

\bpf
Part (ii) is proved in \cite[Lemma 3.4 and Proposition 5.3]{MeSa:11}: see \cite[Proof of Theorem 5.8]{MeSa:11}, and observe that the nontrapping assumption implies that $\alpha$ in \cite{MeSa:11} equals zero. We now prove Part (i).

By the consequence \eqref{eq:Cint2} of the definition of $\Cint$ \eqref{eq:Cint}, 
there exists $v_h \in \cH_h$ such that
\beqs
\N{\cS^*f -v_h}_{\HoDk} \leq \sqrt{2}\Cint h |\cS^*f|_{H^2(\OR)}
\eeqs
(indeed, we can take $v_h = I_h(\cS^* f)$).
By \eqref{eq:S*fkey}, the BVP \eqref{eq:H2equation} is satisfied with $v:= \cS^*f$ and $\widetilde{f}:=f + k^2 n \cS^*f$. Applying the bounds \eqref{eq:Creg} and \eqref{eq:Csol}, we obtain
\begin{align*}
|\cS^*f |_{H^2(\OR)} &\leq \Creg\left( k^2 n_{\max} \N{\cS^*f}_{L^2(\OR)} + \N{f}_{L^2(\OR)}\right.\\ 
&\hspace{3cm}\left.+ \frac{1}{R}\N{\nabla(\cS^*f)}_{L^2(\OR)} + \frac{1}{R^2} \N{\cS^* f}_{L^2(\OR)}\right)\\
&\leq \Creg\Csol kR \left( n_{\max} + \frac{1}{k R\,\Csol} + \frac{1}{kR} + \frac{1}{(kR)^2} \right)\N{f}_{L^2(\OR)},
\end{align*}
and the result \eqref{eq:etabound1} follows from the assumption that $kR\geq k_0 R_0\geq 1$ (see \eqref{eq:assum}).
\epf

\section{Proof of the oscillatory-behaviour bound \eqref{eq:Cosci}}\label{sec:Cosci}

\bth\label{thm:Cosci}
If $\MA, n, $ and $\Omega_-$ are nontrapping (in the sense that the bound \eqref{eq:Csol} holds), 
then the bound \eqref{eq:Cosci} holds, i.e.,
\beq\label{eq:Cosci2}
\vert u\vert_{H^2(\Omega_R)}  \leq \Cosc k\N{u}_{H^1_k(\Omega_R)}.
\eeq
\enth

\ble\label{lem:L2}
To prove Theorem \ref{thm:Cosci},
it is sufficient to prove that there exists $k_0>0$ and $\Cmass=\Cmass(\MA, n, \Omega_-, R)>0$ such that
\beq\label{eq:L2}
\N{u}_{L^2(\Omega_{R+1})} \leq  \Cmass\N{u}_{L^2(\Omega_{R})} \quad \tfa k\geq k_0.
\eeq
\ele

\bpf
We first claim that the map $k\mapsto u$ is continuous from
$(1,\infty)$ to $H^2(\OR)$; indeed, this follows from
the well-posedness of the plane-wave scattering problem of Definition \ref{def:planewave}, $H^2$
regularity, and linearity. Therefore, the function $k \mapsto 
\N{u}_{H^2(\OR)} \big(k \N{u}_{\HoDk}\big)^{-1}$ 
is continuous on $[1,\infty)$, and it is sufficient to prove that the bound \eqref{eq:Cosci2} (i.e., \eqref{eq:Cosci}) holds for $k$ sufficiently large.

Let $\chi\in C^{\infty}(\mathbb{R}^{d})$ be such that $0\leq\chi\leq1$, $\chi =1$ on $\OR$ and $\chi=0$ on $\Rea^d\setminus B_{R+1/2}$.
Applying the $H^2$-regularity results \cite[Theorems 2.3.3.2, 2.4.2.5, and 2.4.2.7]{Gr:85} to $\chi u$ (with these results valid since $\MA$ is Lipschitz, $A_{\min}>0$, both $\Gamma$ and $\GR$ are $C^{1,1}$, and either $\gamma u=0$ or $\partial_{\bn}u=0$ on $\Gamma$),
we obtain, in a similar way to the proof of Theorem \ref{thm:Creg}, that
there exists 
$C_1= C_1(\MA, n, \Omega_-, R)>0$, such that
\beqs
|u|_{H^2(\OR)}\leq C_1 k \N{u}_{H_k^1(\Omega_{R+1})}.
\eeqs
Therefore to prove \eqref{eq:Cosci2} (i.e., \eqref{eq:Cosci}), it is sufficient to prove that there exists 
$C_2= C_2(\MA, n, \Omega_-, R)>0$, such that
\beq\label{eq:osci_2}
\N{u}_{H_k^1(\Omega_{R+1})}\leq C_2 \N{u}_{H_k^1(\Omega_{R})}.
\eeq
We now need to show that we can prove \eqref{eq:osci_2} from \eqref{eq:L2}.
We claim that 
\beq\label{eq:Green}
\N{\nabla u}_{L^2(\Omega_{R+1})} \leq \sqrt{\frac{n_{\max}}{A_{\min}}} k \N{u}_{L^2(\Omega_{R+1})} \quad\tfa k>0.
\eeq
Indeed, applying Green's identity in $\Omega_R$ (which is justified by \cite[Theorem 4.4]{Mc:00} since $u \in H^1(\OR)$) and recalling that either $\gamma u=0$ or $\partial_\bn u=0$ on $\Gamma$, we have that
\beqs
\int_{\Omega_{R+1}} (\MA \nabla u )\cdot\overline{\nabla u} - k^2 n |u|^2 = \Re \int_{\Gamma_{R+1}} \overline{u} \pdiff{u}{r}.
\eeqs
By \eqref{eq:CDtN2}, the right-hand side is $\leq 0$, and \eqref{eq:Green} follows using the inequalities \eqref{eq:Alimits} and \eqref{eq:nlimits}.
Therefore, using \eqref{eq:Green} and \eqref{eq:L2}, 
\beqs
\N{u}_{H_k^1(\Omega_{R+1})}\leq 
\sqrt{\frac{n_{\max}}{A_{\min}}+1}\,\,
k\N{u}_{L^2(\Omega_{R+1})} \leq
\Cmass\sqrt{\frac{n_{\max}}{A_{\min}}+1}\,\,
 k\N{u}_{L^2(\Omega_{R})}
\eeqs
which implies the bound \eqref{eq:osci_2}, and the result follows.
\epf

\subsection{Overview of the ideas used in the rest of this section to prove \eqref{eq:L2}}\label{sec:defect1}
We have therefore reduced proving the oscillatory-behaviour bound \eqref{eq:Cosci}/\eqref{eq:Cosci2} to proving the bound \eqref{eq:L2}, which we prove using \emph{defect measures}. The precise definition of a defect measure is given in Theorem \ref{thm:defect} below, but the idea is that the defect measure of a Helmholtz solution describes 
where the mass of the solution in phase space 
(i.e.~the set of positions $\bx$ and momenta $\bxi$) 
is concentrated in the high-frequency limit.
Two examples of this feature are 

(i) 
the defect measure of the plane wave $u^I(\bx):= \exp (\ri k \bx\cdot\ba)$ is the product of a delta function at $\bxi=\ba$ and Lebesgue measure in $\bx$ (see \eqref{eq:defectpw} below), reflecting the fact that, at high frequency (and in fact at any frequency), all the mass in phase space of the plane wave is travelling in the direction $\ba$,  
and 

(ii) the defect measure of an outgoing solution of the Helmholtz equation is zero on the so-called ``directly incoming set" 
(see Lemma \ref{lem:muI} below), where this set is defined
in \eqref{eq:incomingset} below as points in phase space that don't hit the scatterer when propagated backwards along the flow.

A key feature of the defect measure of a Helmholtz solution is that it is invariant under the Hamiltonian flow defined by the symbol of the PDE, as long as the flow doesn't encounter the scatterer (see Theorem \ref{thm:invariance} below)
This is analogous to results about propagation of singularities of the wave equation, where singularities travel along the trajectories of the flow (the \emph{bicharacteristics}), and the projection of these trajectories in space are the \emph{rays}. 

The main ingredients to our proof of \eqref{eq:L2} are Points (i) and (ii) above, invariance under the flow (away from the scatterer), and then geometric arguments about the rays, using the fact that away from the scatterer the rays are straight lines and the flow has constant speed along the rays (see \eqref{eq:straightlines} below).

To conclude this overview, we direct the reader to  \cite[Chapter 5]{Zw:12} for extensive discussion of defect measures in $\Rea^d$, to  \cite{Bu:02,Mi:00,GaSpWu:20} for material on defect measures on manifolds with boundary, and to  \cite{Bu:97} 
for discussion on the history of defect measures. 

\subsection{Recap of results about defect measures}\label{sec:defect2}

\subsubsection{Symbols and quantisation}

Before defining defect measures, we need to define the functions on
phase space (i.e.~the set of positions $\bx$ and momenta $\bxi$) that
the defect measure can act upon by dual pairing. These
functions are called \emph{symbols}, defined as functions on the
\emph{cotangent bundle} $T^*\Omega_+$. 
Recall the definition of the cotangent bundle of
  $\Rea^d$: 
  \beqs
T^*\Rea^d :=\Rea^d \times (\Rea^d)^*;
\eeqs
for our purposes, we can consider $T^*\Rea^d$ as $\{(\bx,\bxi) : \bx\in \Rea^d, \bxi\in\Rea^d\}$, i.e.~the set of positions $\bx$ and momenta $\bxi$.
On $T^*\Rea^d$, the \emph{quantisation} of a symbol $b(\bx,\bxi) \in C_{\rm comp}^\infty(T^*\Rea^d)$ is defined by
\beq\label{eq:quant}
b\big(\bx, (\ri k)^{-1}\partial_\bx\big)u(\bx):= \frac{k^d}{(2\pi)^d} \int_{\Rea^d}\int_{\Rea^d} \re^{\ri k (\bx-\by)\cdot \bxi} \,b(\bx,\bxi) u(\by) \, \rd \by \,\rd \bxi;
\eeq
see, e.g., \cite[\S4]{Zw:12}. The same definition holds for symbols supported away from the boundary of $\overline{\Omega}_+$.
We omit the analogous definition near the boundary since it is more involved; see \cite[\S4.2]{Bu:02}
(where it involves the so-called \emph{compressed cotangent bundle} of $\Omega_+$,
 $T^*_{\rm b}\overline{\Omega_+}$) and \cite[\S1.2]{Mi:00}. We will not, in any event,
require any specifics of the measure at the boundary in proving
Theorem \ref{thm:Cosci}.
 
 \subsubsection{Existence of defect measures}

\bth\mythmname{{\bf Existence of defect measures \cite[Theorem 5.2]{Zw:12}, \cite[\S4.2]{Bu:02}}}\label{thm:defect}
Suppose $\{v(k)\}_{k_0\leq k<\infty}$ is a collection of functions that is uniformly locally bounded in $L^2(\Omega_+)$, i.e.~given $\chi \in C_{\rm comp}^\infty(\Rea^d)$ there exists $C>0$, depending on $\chi$ and $k_0$ but independent of $k$, such that 
\beq\label{eq:defectv}
\N{\chi v(k)}_{L^2(\Omega_+)}\leq C \quad\tfa k \geq k_0.
\eeq
Then there exists a sequence $k_{\ell}\rightarrow\infty$ and a non-negative Radon measure $\mu$ on
$T^*_{\rm b}\overline{\Omega_+}$ (depending on $k_{\ell}$) such that, for any symbol $b(\bx,\bxi)\in C_{\rm comp}^\infty(T^*_{\rm b}\overline{\Omega_+})$
\beq\label{eq:defect}
\big\langle b\big(\bx,(\ri k_\ell)^{-1}\partial_{\bx}\big)v(k_\ell),v(k_{\ell}))\big\rangle_{\Omega_+}\longrightarrow\int b\ \rd\mu \quad \tas \ell \rightarrow \infty.
\eeq
\enth
In the case of a plane wave $u^I(\bx):= \exp (\ri k\bx\cdot\ba)$ with $|\ba|=1$, a direct calculation using \eqref{eq:quant} and the definition of the Fourier transform shows that, for all $k$,
\begin{align}\nonumber
\big\langle b\, u^I, u^I\big\rangle_{\Rea^d} &:= \frac{k^d}{(2\pi)^d} \int_{\Rea^d}\int_{\Rea^d}\int_{\Rea^d}
\,\re^{\ri k(\bx-\by)\cdot\bxi}\, \re^{\ri k \by\cdot\ba} \,\re^{-\ri k \bx\cdot\ba}b(\bx,\bxi)\,\rd \bxi\, \rd \by \,\rd \bx
\\
&=\int_{\Rea^d} b(\bx,\ba)\,\rd \bx;\label{eq:defectpw}
\end{align}
i.e.~for any sequence $k_{\ell}\rightarrow\infty$, the corresponding defect measure of $u^I$ is the product of the Lebesgue measure in $\bx$ by a delta measure at $\bxi=\ba$; we therefore talk about \emph{the} (as opposed to \emph{a}) defect measure of $u^I$.

The next lemma proves that, if $u$ is the solution of the plane-wave scattering problem
and $\chi$ is an arbitrary cut-off function, then $\chi u$ is uniformly bounded in $k$ (on compact subsets of $\Omega_+$); existence of a defect measure of $u$ then follows from Theorem \ref{thm:defect}.
In the rest of this section, to emphasise the $k$-dependence of $u$, we write $u=u(k)$.

\ble\label{lem:uniformbound}
Let $u(k)$ be the solution of the plane-wave scattering problem of Definition \ref{def:planewave}.
Assume that $\MA, n, $ and $\Omega_-$ are 
nontrapping. Then there exists $C(\MA,n,\Omega_-, R, k_0)>0$ such that
\beq\label{eq:uniformbound}
\N{u(k)}_{L^2(\OR)} \leq C
 \quad\tfa k \geq k_0.
\eeq
\ele

\bpf
Let $\chi\in C^\infty_{\rm comp}(\Rea^d)$ be such that $\chi=1$ in a neighbourhood of the scatterer $\Omegascat$.
Let $v:= u^S + \chi u^I$, so that $u= (1-\chi)u^I + v$. 
Since $\|u^I(k)\|_{L^2(\OR)}\leq C_1(R)$ for all $k>0$, the result \eqref{eq:uniformbound} will follow if we prove a uniform bound on $\|v(k)\|_{L^2(\OR)}$.
The definition of $v$ implies that $v$ satisfies the Sommerfeld radiation condition, either $\gamma v=0$ or $\partial_\bn v=0$ on $\Gamma$, and, with $\cL_{\MA,n}w:= \nabla\cdot(\MA \nabla w ) + k^2 n w$ and $[A,B]:=AB-BA$,
\beqs
\cL_{\MA,n}v = - \cL_{\MA,n}\big((1-\chi)u^I\big) = \big[ \cL_{\MA,n}, \chi \big]u^I - (1-\chi) \cL_{\MA,n} u^I = \big[ \cL_{\MA,n}, \chi \big]u^I,
\eeqs
since $\cL_{\MA,n}u^I=0$ when $1-\chi \neq 0$. By explicit calculation, using the fact that $u^I(\bx)= \exp (\ri k \bx\cdot\ba)$,
\beqs
\N{ \big[ \cL_{\MA,n}, \chi \big]u^I }_{L^2(\OR)} \leq C_1 
\, k,
\eeqs
where $C_1$ depends on $\|\MA\|_{L^\infty(\OR)}, \|\nabla\MA\|_{L^\infty(\OR)},$ and $\chi$, but is independent of $k$.
The nontrapping bound \eqref{eq:Csol} then implies that $\|v(k)\|_{L^2(\OR)}\leq C_2$ with $C_2$ independent of $k$, and the result follows.
\epf

 \subsubsection{Support and invariance properties of defect measures}

Recall that the semi-classical principal symbol of the Helmholtz equation \eqref{eq:PDE} is given by
\beq\label{eq:principal_symbol}
p(\bx,\bxi):= \sum_{i=1}^d\sum_{j=1}^{d} A_{ij}(\bx)\xi_i \xi_j - n(\bx)
\eeq
(see, e.g., \cite[Page 281]{Zw:12}).
In our arguments below we only consider points $(\bx,\bxi)$ in phase space when $p=0$; this is because of the following result.

\bth\mythmname{{\bf Support of defect measure \cite[Theorem 5.4]{Zw:12}, \cite[Equation 3.17]{Bu:02}}}
 \label{thm:support}
Suppose $u(k)$ satisfies \eqref{eq:uniformbound}, and let $\mu$ be any defect measure of $u(k)$.
Then $\supp \mu \subset \{ (\bx,\bxi) : p(\bx,\bxi)= 0\}.$
\enth

As an illustration of this, the plane wave $u^I(\bx):= \exp (\ri k\bx\cdot\ba)$ with $|\ba|=1$ is solution of the Helmholtz equation \eqref{eq:PDE} with $\MA=\MI$ and $n=1$, and hence $p=|\bxi|^2-1$ in this case. 
By \eqref{eq:defectpw}, the defect measure of $u^I$ is the product of Lebesgue measure in $\bx$ and a delta function at $\bxi=\ba$, and thus is supported in $|\bxi|=1$, i.e., $p=0$, as expected from Theorem \ref{thm:support}.

The final result about defect measures that we need is their invariance under the flow (away from the scatterer). This result is Theorem \ref{thm:invariance} below; to state it, we first need to define the flow.

Away from $\Gamma$, and provided that $A$ and $n$ are both $C^{1,1}$, the flow $\varphi_t$ is defined as follows: given $\rho= (\bx_0,\xi_0)$, $\varphi_t(\rho):= (\bx(t),\bxi(t))$ where $(\bx(t),\bxi(t))$ is the solution of 
the Hamiltonian system
\beq\label{eq:Hamilton}
\dot{x_i}(t) = \partial_{\xi_i}p\big(\bx(t), \bxi(t) \big), \qquad
\dot{\xi_i}(t)
 = -\partial_{x_i}p\big(\bx(t), \bxi(t) \big),
\eeq
with initial condition $(\bx(0), \bxi(0))= (\bx_0, \bxi_0)$,
where the Hamiltonian equals $p$ defined by \eqref{eq:principal_symbol}.
Near both $\Gamma$ and places where $A$ and $n$ are not $C^{1,1}$, the definition of $\varphi_t$ is more involved -- this is to account for reflection or refraction. However, we do not need this definition in what follows, since our arguments take place away from these regions. In fact our arguments take place away from the scatterer $\Omegascat$. Outside $\Omegascat$, $\MA=\MI$, and $n=1$; thus $p(\bx,\bxi)= |\bxi|^2-1$. From \eqref{eq:Hamilton}, the flow satisfies
$\dot{x_i} =2 \xi_i$ and $\dot{\xi_i}=0$
and is therefore given by the straight-line motion
\beq\label{eq:straightlines}
\bx=\bx_0+2t\bxi_0, \quad \bxi=\bxi_0.
\eeq
The arguments below consider the flow with speed $2$ (i.e.~with $|\bxi_0|=1$). This is without loss of generality, since 
away from $\Omegascat$  Theorem \ref{thm:support} implies
  that $\mu$ is only non-zero when
$|\bxi|=1$. 

Both in the next result and later, we let $\pi_\bx$ denote projection in the $\bx$ variables, i.e.~$\pi_\bx((\bx,\bxi))=\bx$.

\bth[Invariance of defect measure under the flow away from the scatterer]
 \label{thm:invariance}
 Suppose that $u(k)$ satisfies \eqref{eq:uniformbound}, and let $\mu$ be any defect measure of $u(k)$.
 If $A\subset T^* \Rea^d$ is such that $\pi_\bx (\varphi_s(A)) \cap
\Omegascat= \emptyset$ for $s$ between $0$ and $t$, 
(i.e.~the flow acting on $A$ doesn't hit the scatterer from time $0$ to time $t$), then
\beq\label{eq:invariance}
\mu( \varphi_t(A))= \mu(A).
\eeq
\enth
\bpf
In the absence of the scatterer, invariance of the measure under the flow is the statement that, for   $b\in C_{\rm comp}^\infty(T^*\Rea^d)$, 
\beq\label{eq:Einvariance}
\partial_s \bigg( \int (b\circ \varphi_{-s})(\rho)\, \rd \mu\bigg) = 0 \quad \tfa s,
\eeq
and this is proved in \cite[Theorem 5.4]{Zw:12}, \cite[Proposition 4.4]{Bu:02}. 
For this result to hold in the presence of the scatterer in a time interval $0\leq s\leq t$, we need the spatial projection of the integrand in \eqref{eq:Einvariance} to not be supported during this time interval 
 on $\Omegascat$, i.e., we need the condition that
\beq\label{eq:Econdition}
\pi_\bx\big(\supp(b\circ \varphi_{-s})\big) \cap \Omegascat = \emptyset \quad\tfor 0\leq s\leq t.
\eeq
Under this condition, \eqref{eq:Einvariance} implies that
\beq\label{eq:Estar}
 \int b(\rho)\, \rd \mu= \int (b\circ \varphi_{-s})(\rho)\, \rd \mu \quad\tfa 0\leq s\leq t.
\eeq

Let $1_A$ denote the indicator function of a set $A$. By approximating $1_A$ by smooth symbols, \eqref{eq:Estar} holds with $b(\rho)=1_A(\rho)$, provided that the condition \eqref{eq:Econdition} holds.
Since $\varphi_{-s}(\rho) \in A$ iff $\rho \in \varphi_s(A)$, we have
\beqs
\pi_\bx\big(\supp(1_A\circ \varphi_{-s})\big)= \pi_\bx\big(\supp(1_{\varphi_{s}(A)})\big) = \pi_\bx\big(\varphi_s(A)\big),
\eeqs
and thus \eqref{eq:Econdition} holds by the assumption in the statement of the theorem.

Therefore, \eqref{eq:Estar} implies that, for all $0\leq s\leq t$,
\begin{align*}
\int 1_A(\rho)\, \rd \mu &= \int 1_A(\varphi_{-s}(\rho))\, \rd \mu = \int 1_{\varphi_s(A)}(\rho)\, \rd \mu,
\end{align*}
i.e.
\beqs
\mu(A) = \mu\big(\varphi_s(A)\big) \quad\tfa 0\leq s\leq t,
\eeqs
which implies \eqref{eq:invariance}.
\epf  
 
\subsection{Proof of \eqref{eq:L2} using defect measures}

The following lemma reduces proving the bound \eqref{eq:L2} to proving a statement about defect measures. 

\ble\label{lem:Cosci}
Let $0<R_0<R$ be such that $\Omegascat \subset \subset B_{R_0}$.
If every defect measure of $u$ is non-zero and there exists $C_{R,R_0}>0$ such that, for every defect measure $\mu$ of $u$,
\beq \label{eq:reduc_defect}
\mu(T^*\Omega_{R+2}) \leq C_{R,R_0} \mu(T^*\Omega_{R_0}),
\eeq
then the bound \eqref{eq:L2} holds.
\ele

\bpf
We prove the contrapositive.  Suppose \eqref{eq:L2} fails; we aim
  to exhibit a defect measure associated to $u$ for which \eqref{eq:reduc_defect} fails.
Then, for any $C_1>0$, there exists a sequence $(k_n)_{n=1}^\infty$, with $k_n \rightarrow \infty$, such that
\beq\label{eq:DM1}
\Vert	u(k_n) \Vert_{L^2(\Omega_{R+1})} \geq C_1  \Vert u(k_n) \Vert_{L^2(\Omega_{R})};
\eeq
we choose $C_1:= 2C_{R,R_0}$.
By Lemma \ref{lem:uniformbound}, the sequence $\{u(k_n)\}_{n=1}^\infty$ is locally uniformly bounded and
Theorem \ref{thm:defect} implies
  that, by passing to a subsequence, there exists a defect measure
  $\mu$ of $u$ associated to the subsequence, which we again denote $k_n$.
Let $\chi_0, \chi_1 \in C^\infty(\Rea^d)$  be such that  $0\leq \chi_{0}, \chi_1 \leq 1$, and 
$$
\supp \chi_1 \subset B_{R+2}, \quad\chi_1 = 1 \text{ in } B_{R+1},  \quad \supp \chi_0 \subset B_{R},\quad \chi_0 = 1 \text{ in } B_{R_0}.
$$
The bound \eqref{eq:DM1} then implies that 
\beq\label{eq:Euan2}
\Vert	\chi_1 u(k_n) \Vert_{L^2(\Omega_+)} \geq  2C_{R,R_0} \Vert \chi_0  u(k_n) \Vert_{L^2(\Omega_+)}.
\eeq
Passing to the limit $n\rightarrow \infty$ and using the property of defect measure \eqref{eq:defect}, we obtain that
\beqs
\int \chi_1^2 \,\rd \mu \geq 2C_{R,R_0} \int \chi_0^2 \,\rd \mu.
\eeqs
The definitions of $\chi_{0}$ and $\chi_{1}$ imply that
\beqs
\int \chi_0^2 \, \rd \mu \geq \int 1_{T^*\Omega_{R_0}}\, \rd \mu = \mu( T^*\Omega_{R_0})
\eeqs
and 
\beqs
\int \chi_1^2 \, \rd \mu \leq \int 1_{T^*\Omega_{R+2}}\, \rd \mu = \mu( T^*\Omega_{R+2});
\eeqs
hence
$$
\mu(T^*\Omega_{R+2}) \geq 2C_{R,R_0} \mu(T^*\Omega_{R_0}),
$$
contradicting \eqref{eq:reduc_defect}.
\epf

Before using Lemma \ref{lem:Cosci} to prove \eqref{eq:L2}, we prove a result (Lemma \ref{lem:muI} below) about the structure of $\mu$, exploiting the fact that $u=u^I + u^S$ with $u^S$ is outgoing (in the sense that it satisfies the Sommerfeld radiation condition \eqref{eq:src}). To make use of this outgoing property, we need to define appropriate notions of incoming and outgoing for elements of phase space.
Let $\mathcal I$ denote the \emph{directly incoming set} defined by
\beq\label{eq:incomingset}
\cI:= \bigg\{\rho\in T^{*}(\Omega_+{\setminus\Omegascat}),\text{ s.t. }\pi_\bx\bigg(\bigcup_{t\geq0}\varphi_{-t}(\rho)\bigg)\cap\Omegascat=\emptyset\bigg\};
\eeq
where recall that $\pi_\bx$ denotes projection in the $\bx$ variables. 
That is, $\cI$ is everything that never hits the scatterer under backward flow. 
Let
\beqs
\Gamma_+ := (T^* \Omega_+) \backslash \cI.
\eeqs
 These definitions of $\cI$ and $\Gamma_+$ do not require the generalized bicharacteristic flow $\varphi_t$ to be defined in $T^* \Omegascat$, but when the flow is defined everywhere, $\Gamma_+$ is the forward generalized bicharacteristic flowout of $\Omegascat$, that is
\beqs
\Gamma_+= \bigg\{ \bigcup_{t\geq 0} \varphi_t (\rho) \,\, : \,\, \rho \in T^* \Omegascat\bigg\} \text{ when }\varphi_t\text{ is defined everywhere.}
\eeqs

The following lemma uses outgoingness of $u^S$ to show that, given a set $E$ in phase space, the mass of $u$ lying over $E$ is \emph{either} in the forward flowout $\Gamma_+$ \emph{or} associated to the incident wave $u^I$.

\ble \label{lem:muI}
For any Borel set $E\subset T^* \Omega$,  $\mu(E \setminus \Gamma_+) = \mu^I (E \setminus \Gamma_+) $, where $\mu$ is any defect measure of $u$, and $\mu^I$ is the defect measure of $u^I$.
\ele
\bpf
Let $k_\ell$ be the sequence associated to the particular defect measure of $u$. By Lemma \ref{lem:uniformbound}, $u^S(k_\ell)$ is uniformly locally bounded, and so there exists a subsequence $k_{\ell_m}$ and a defect measure associated to $u^S$, denoted by $\mu^S$. Then, by linearity and \eqref{eq:defect}, $\mu= \mu^S+ \mu^I$.
It is therefore sufficient to prove that $\mu^S(E\setminus\Gamma_+)=0$. But, by the definition of $\Gamma_+$, $E\setminus\Gamma_+ \subset \cI$, and $\mu^S(\cI)=0$ by \cite[Proposition 3.5]{Bu:02}, \cite[Lemma 3.4]{GaSpWu:20}, since $u^S$ is outgoing.
\epf

\bpf[Proof of Theorem \ref{thm:Cosci}]
By Lemmas \ref{lem:L2} and \ref{lem:Cosci} it is sufficient to prove the bound \eqref{eq:reduc_defect}
(observe that the hypothesis in Lemma \ref{lem:Cosci} that every defect measure of $u$ is non-zero holds by Lemma \ref{lem:muI} since $\mu^I(\cI)\neq 0$).
Let $\Rscat:= \max_{\bx\in \Omegascat}|\bx|$. We claim that it is sufficient to show that, for any $\rho > \Rscat$ there exists $\eps=\eps(\Rscat,\rho)$ , with $\eps(\Rscat,\rho)$ is an increasing function of $\rho$, 
and $C=C(\rho,\eps)>0$ such that
\beq\label{eq:DM2}
\mu(T^*(B_{\rho+\eps} \setminus B_{\rho})) \leq C(\rho,\eps) \mu(T^*\Omega_\rho).
\eeq
Indeed, we now show that the bound \eqref{eq:reduc_defect} then follows by using \eqref{eq:DM2} repeatedly. Since $\eps(\Rscat,\rho)$ is an increasing function of $\rho$, if
$\eps^* := \eps(\Rscat, R_0)$, then \eqref{eq:DM2} implies, with $C(\rho):= C(\rho, \eps(\Rscat,\rho))$, 
\beq\label{eq:DM2_upgr}
\mu(T^*(B_{\rho+\eps^*} \setminus B_{\rho})) \leq C(\rho)\, \mu(T^*\Omega_\rho) \quad\tfa \rho \geq R_0.
\eeq
The bound  \eqref{eq:reduc_defect} then follows by applying \eqref{eq:DM2_upgr} with $\rho=R_0$, $\rho = R_0+ \eps^*$, \ldots ,  $\rho= R_0 + m\eps^*$, where $m = \lceil (R+2-R_0)/\eps^* \rceil$.

It is therefore sufficient to prove the bound \eqref{eq:DM2}; we introduce the notation that $A:=B_{\rho+\eps} \setminus B_{\rho}$, and observe that \eqref{eq:DM2} then reads $\mu(T^*A) \leq C(\rho,\eps) \mu(T^*\Omega_\rho)$.
We prove this bound by combining the following three inequalities:
\beq\label{eq:DMkey1}
\mu(T^*A) \leq \mu(T^*A \cap \Gamma_+) + \mu_I(T^*A) = \mu(T^*A \cap \Gamma_+)+ |A|
\eeq
 (where $|\cdot|$ denotes Lebesgue measure in $\Rea^d$),
 \beq\label{eq:DMkey2} 
\mu(T^*A \cap\Gamma_+)\leq \mu (T^*(B_{\rho} \setminus B_{\rho_0}))\leq \mu (T^*\Omega_\rho),
\eeq  
where $\rho_0:= (\rho+\Rscat)/2$, 
and 
\beq\label{eq:DMkey3}
\mu(T^*\Omega_\rho) \geq \delta |\Omega_\rho|
\eeq
for some $\delta>0$. Indeed, using \eqref{eq:DMkey1}, \eqref{eq:DMkey2}, and \eqref{eq:DMkey3}, we have
\beqs
\mu(T^*A) \leq \Big(1+ \abs{A}(\delta \abs{\Omega_\rho})^{-1}\Big)
\mu(T^*\Omega_\rho), 
\eeqs
which is \eqref{eq:DM2}. We prove \eqref{eq:DMkey1} and \eqref{eq:DMkey3} using Lemma \ref{lem:muI} and the structure of $\mu^I$, and \eqref{eq:DMkey2} using invariance of defect measures under the flow outside of $T^* \Omegascat$ (i.e.~Theorem \ref{thm:invariance}).

\paragraph{Proof of \eqref{eq:DMkey1}}
Lemma \ref{lem:muI} implies that
\beqs
\mu(T^*A) = \mu(T^*A\cap \Gamma_+) + \mu(T^*A\setminus \Gamma_+) \leq \mu(T^*A \cap \Gamma_+) + \mu_I(T^*A).
\eeqs
By \eqref{eq:defectpw}, $\mu_I$ is a $\delta$-measure on $\bxi=\ba$ times Lebesgue measure in $\bx,$ so $\mu_I(T^*A) = \abs{A},$ 
(where $|\cdot|$ denotes Lebesgue measure in $\Rea^d$) and \eqref{eq:DMkey1} follows.

\paragraph{Proof of \eqref{eq:DMkey2}}
Recall that, for $X\subset\subset \Rea^d\setminus\overline{\Omegascat}$,
\beqs
S^*X:= \big\{ (\bx,\bxi) : \bx\in X,\, \bxi \in \Rea^d \text{ with} \,|\bxi|=1\big\},
\eeqs
and observe that, by Theorem \ref{thm:support}, $\mu(T^*A\cap \Gamma_+)= \mu(S^*A\cap \Gamma_+)$ and $\mu (T^*(B_{\rho} \setminus B_{\rho_0})) = \mu (S^*(B_{\rho} \setminus B_{\rho_0}))$; we therefore only need to prove that
 \beq\label{eq:DMkey2bis} 
\mu(S^*A \cap\Gamma_+)\leq \mu (S^*(B_{\rho} \setminus B_{\rho_0})).
\eeq 
We first introduce some notation that allows us to bound
  $\mu(S^*A \cap\Gamma_+)$ using only the invariance of defect measure \eqref{eq:invariance}
in the exterior of $\Omegascat$. 
Given $\bbb \in \mathbb R ^d$ with $|\bbb|=1$
and $\widetilde{\rho}>\Rscat$, let $\Omegascatb\subset \Rea^d$
and $\Lambda_{\text{sc},\widetilde{\rho},\bbb} \subset S^* \Omega_+$
be defined by \beqs \Omegascatb:= \Big( \bigcup_{t \geq 0}
\big(\Omegascat + t \bbb\big) \Big)\cap \Omega_{\widetilde{\rho}}
\quad\tand \quad \Lambda_{\text{sc},\rho+\eps,\bbb} := \Omegascatb
\times \{\bbb\}; \eeqs i.e.~$\Omegascatb$ equals the union of all
possible translations of $\Omegascat$ in the direction $\bbb$,
intersected with $\Omega_{\widetilde{\rho}}$, and
$\Lambda_{\text{sc},\widetilde{\rho},\bbb}$ equals these points paired
with the direction $\bbb$.
By \eqref{eq:straightlines}, the spatial projections of the flow
outside $\Omegascat$ are straight lines, and thus
\beqs \Gamma_+ \cap S^*
\Omega_{\widetilde{\rho}} \cap \{ \bxi = \bbb \} = \Big\{ (\bx, \bbb)
\in S^* \Omega_{\widetilde{\rho}} : \exists s \geq 0 \text{ s.t. }
\bx-s\bbb \in \Omegascat \Big\}.  
\eeqs 
Therefore
\beq \label{eq:cover1} \Gamma_+ \cap S^* \Omega_{\widetilde{\rho}}
\cap \{ \bxi = \bbb \} \subset
\Lambda_{\text{sc},\widetilde{\rho},\bbb}, \qquad \Gamma_+ \cap S^*
\Omega_{\widetilde{\rho}} \subset \bigcup_{\bbb \in \mathbb R ^ d,  |\bbb|=1}
\Lambda_{\text{sc},{\widetilde{\rho}},\bbb}, 
\eeq 
  and thus, for any $\eps>0$,
  \beq\label{eq:flowout} 
  S^*A \cap \Gamma_+ = S^*A \cap
  S^*\Omega_{\rho+\eps}\cap \Gamma_+\subset S^*A \cap \bigg(
  \bigcup_{\bbb \in \mathbb R ^ d, |\bbb|=1}
  \Lambda_{\text{sc},\rho+\eps,\bbb}\bigg),  
  \eeq
  Recall
  that $\rho_0:= (\rho+ \Rscat)/2$. Let \beq\label{eq:t0} t_0:=
  \frac{\rho_0-\Rscat}{4} = \frac{\rho-\Rscat}{8} \eeq and
  \beq\label{eq:eps} \eps := -\rho + \sqrt{ \Rscat^2+\left(
      \frac{\rho-\Rscat}{4}+ \sqrt{\rho^2-\Rscat^2}\right)^2}; \eeq
  observe that $\eps>0$ and $\eps$ is an increasing function of
  $\rho$, as claimed underneath \eqref{eq:DM2}.  We now claim that,
  with these definitions of $t_0$ and $\eps$,
  \beq \label{eq:claim1} \bigcup_{0\leq t \leq t_0} \varphi_t
  \big(S^*(B_\rho \setminus B_{\rho_0})\big)
 \cap \Omegascat = \emptyset
 \eeq
(i.e., the forward flowout of the annulus $B_\rho\setminus B_{\rho_0}$ does not hit the scatterer for $0\leq t\leq t_0$) and
\beq\label{eq:claim2} S^*A\cap\bigg( \bigcup_{\bbb \in \mathbb R ^ d, |\bbb|=1}
\Lambda_{\text{sc},\rho+\eps,\bbb} \bigg) \subset
\varphi_{t_0}\big(S^*(B_\rho \setminus B_{\rho_0})\big). \eeq
 (Since $S^*A \cap \Gamma_+$ is
contained in the left-hand side of \eqref{eq:claim2} by
\eqref{eq:flowout}, \eqref{eq:claim2} says that the forward flowout
of $B_\rho \setminus B_{\rho_0}$ in time $t_0$ covers all points in
$S^*A$ that are ever reached by flowout from $T^* \Omega_\scat$.)
Outside $\Omegascat$ the flow has speed 2 and its spatial projections
are straight lines. Therefore \eqref{eq:claim1} is ensured if
$t_0< (\rho_0-\Rscat)/2$, which is ensured by \eqref{eq:t0}.

\begin{figure}
\centering{
\includegraphics[width=0.8\textwidth]{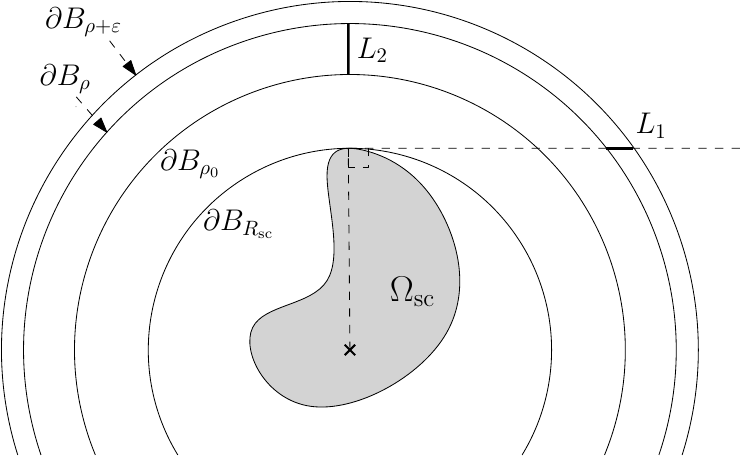}
}
\caption{Figure showing the lengths $L_1$ and $L_2$ defined by \eqref{eq:lengths2}.}\label{fig:1}
\end{figure}

We now show that  \eqref{eq:claim2} holds. 
Since
$$
(\bx, \bbb) = (\bx-2t_0 \bbb + 2t_0\bbb, \bbb) = \varphi_{t_0}(\bx-2t_0 \bbb, \bbb), $$
\eqref{eq:claim2} follows from showing that $(\bx-2t_0 \bbb, \bbb) \in S^*(B_\rho \setminus B_{\rho_0})$, 
i.e.~$\bx-2t_0 \bbb \in B_\rho \setminus B_{\rho_0}$, for all $(\bx, \bbb)$ belonging to the left-hand side of \eqref{eq:claim2}. For such $(\bx, \bbb)$, by definition, 
\beq \label{eq:crown_out}
\rho \leq |\bx| \leq \rho + \epsilon, \text{ and }\bx - s \bbb\in \Omegascat
\eeq
for some $s \geq 0$.
We now claim that for such $(\bx,\bbb)$, 
\beqs
\bx - \ell \bbb \in B_\rho \setminus B_{\rho_0} \quad\tfa \quad L_1 < \ell \leq L_2, 
\eeqs
where 
\beq\label{eq:lengths2}
L_1 :=  \sqrt{(\rho+\eps)^2-\Rscat^2} - \sqrt{\rho^2-\Rscat^2}, \qquad L_2 := \rho - \rho_0.
\eeq
This is because, on the one hand, a ray of length $>L_1$ starting from a point $\bx$ in a direction $-\bbb$, with $(\bx, \bbb)$ satisfying \eqref{eq:crown_out}, will automatically
enter $B_\rho$. Indeed, the longest such ray that does not intersect $B_\rho$ has length $L_1$, as shown in Figure \ref{fig:1}. 
On the other hand, a ray of length $\leq L_2$
starting from a point $\bx$ in a direction $-\bbb$, with $(\bx, \bbb)$ satisfying \eqref{eq:crown_out}, will not intersect $B_{\rho_0}$. Indeed,
the shortest such ray that enters $\overline{B_{\rho_0}}$ has length $L_2$, as shown in Figure \ref{fig:1}.
It is then straightforward to check that $L_1 < 2t_0 \leq L_2$ when $t_0$ is given by \eqref{eq:t0} and $\eps$ is given by \eqref{eq:eps}, so that \eqref{eq:claim2} holds.

We now prove the bound \eqref{eq:DMkey2bis} on $\mu(S^*A \cap\Gamma_+)$ using \eqref{eq:claim1} and \eqref{eq:claim2}.
Because of \eqref{eq:claim1}, we can use \eqref{eq:invariance} to find that 
\beqs
\mu \big( \varphi_{ t_0} (S^*(B_{\rho} \setminus B_{\rho_0}) \big) = \mu (S^*(B_{\rho} \setminus B_{\rho_0}));
\eeqs
using this with \eqref{eq:flowout} and \eqref{eq:claim2}, 
we obtain \eqref{eq:DMkey2bis}, and thus \eqref{eq:DMkey2}.

\paragraph{Proof of \eqref{eq:DMkey3}}
Using Lemma \ref{lem:muI} and the structure of $\mu_I$, we have
\begin{align}
\mu(T^*\Omega_\rho) &\geq \mu(T^*\Omega_\rho \setminus \Gamma_+) = \mu_I(T^*\Omega_\rho \setminus \Gamma_+) \nonumber
\\ \nonumber
&
=\mu_I\big((T^*\Omega_\rho \setminus \Gamma_+)\cap\{ \bxi=\ba\}\big) + \mu_I\big((T^*\Omega_\rho \setminus \Gamma_+)\cap\{ \bxi \neq \ba\}\big)\\
& =\Big|  \pi_\bx \Big( (T^*\Omega_\rho \setminus \Gamma_+)\cap\{ \bxi=\ba\}
 \Big)\Big|. \label{eq:DM5.0}
\end{align}
Since
\beqs
\pi_\bx \Big( (T^*\Omega_\rho \setminus \Gamma_+)\cap\{ \bxi=\ba\}  \Big) \cup \pi_\bx \Big(( T^*\Omega_\rho \cap \Gamma_+)\cap\{ \bxi=\ba\}  \Big) \supset \Omega_\rho.
\eeqs
we obtain
\beq\label{eq:DM5}
\Big|  \pi_\bx \Big( (T^*\Omega_\rho \setminus \Gamma_+)\cap \{ \bxi=\ba\}  \Big)\Big|\geq |\Omega_\rho| - \Big|  \pi_\bx \Big(( T^*\Omega_\rho \cap \Gamma_+)\cap\{\bxi=\ba \} \Big)\Big|.
\eeq
By the first inclusion in \eqref{eq:cover1},
\beq\label{eq:DM6}
\big|\pi_\bx \big( (T^*\Omega_\rho \cap \Gamma_+)\cap\{ \bxi=\ba \} \big) \big| \leq \big|\OmegascataR\big|,
\eeq
with this inequality expressing the fact that any parts of the scattered wave travelling in direction $\ba$ must lie in $\OmegascataR$.
Combining \eqref{eq:DM5} with \eqref{eq:DM6} yields
\beq\label{eq:DM6a}
\Big|  \pi_\bx \big( (T^*\Omega_\rho \setminus \Gamma_+)\cap \{ \bxi=\ba\}  \big)\Big|\geq |\Omega_\rho| - |\OmegascataR|.
\eeq
Since $ \OmegascataR \subsetneq \Omega_\rho$, there exists $\delta>0$ such that $|\Omega_\rho| - |\OmegascataR| \geq \delta |\Omega_\rho|$, and thus \eqref{eq:DM5.0} and \eqref{eq:DM6a} imply that \eqref{eq:DMkey3} holds; the proof is complete.
\epf

\bre[What if impedance boundary conditions are imposed on $\GR$?]
\label{rem:impedance}
If the impedance boundary condition $\partial_\bn u^S - \ri k u^S=0$ is imposed on $\GR$ (as an approximation of $\DtN$), then there are additional reflections on $\GR$ \cite{Mi:00}, \cite[\S2]{GaLaSp:21}
$\mu^S$ has support on the incoming set, and Lemma \ref{lem:muI} no longer holds.
\ere
\bre[Proving Theorem \ref{thm:Cosci} in the trapping case]
\label{rem:trapping2}
In the trapping case, $\|u(k)\|_{L^2(\OR)}$ may no longer be uniformly
bounded, as it is in Lemma \ref{lem:uniformbound}, since
\eqref{eq:Csol} no longer holds with $\Csol$ bounded independently of
$k$.  If a subsequence of $k$'s exists along which
$\|u(k)\|_{L^2(\OR)}$ is uniformly bounded, we may obtain a
contradiction by the same argument as above by considering this
subsequence.  Thus, we can assume, without loss of generality,
that $\|u(k)\|_{L^2(\OR)} \to\infty.$
Now instead of defining defect measures of $u(k)$, one can instead define defect measures of $u(k)/\|u(k)\|_{L^2(\OR)}$.
If $R$ is sufficiently large, then the bound in \cite[Theorem 1.1]{CaVo:02} 
(i.e.~the fact that the nontrapping cut-off resolvent estimate holds, even under trapping, if the supports of the cut-offs on both sides are sufficiently far away from the scatterer)
implies that $v(k):= u(k)/\|u(k)\|_{L^2(\OR)}$ satisfies \eqref{eq:defectv}. 
Any defect measure of $v(k)$ is then immediately non-zero, since $\mu(\chi^2) \geq 1$ for any $\chi$ with $\supp \chi \supset B_R$.
Lemma \ref{lem:Cosci} goes through as before after multiplying both sides of \eqref{eq:Euan2} by $\|u(k)\|_{L^2(\OR)}^{-2}$.
The main change needed to the rest of the proof is to take into account the fact that a defect measure of $u^I(k)/\|u(k)\|_{L^2(\OR)}$ 
is zero when $\|u(k)\|_{L^2(\OR)}$ grows through the sequence $k_\ell$ associated with that measure.
In this situation, however, the bound \eqref{eq:DMkey1} becomes
$\mu(T^* A) \leq \mu(T^*A \cap \Gamma_+)$; combining this with \eqref{eq:DMkey2} we obtain $\mu(T^* A) \leq 2\mu(T^*\Omega_R)$, from which the key bound \eqref{eq:DM2} (and hence the result of the theorem) follows.
\ere
  
\section{Proof of Theorems \texorpdfstring{\ref{thm:main1}}{} and  \texorpdfstring{\ref{thm:main2}}{}}

\ble[Aubin-Nitsche analogue via elliptic projection]\label{lem:AN}
Assuming that the Galerkin solution $u_h$ to the variational problem \eqref{eq:FEM} exists, if 
\beq\label{eq:C1}
hk^2 \eta(\cH_h) \leq \cC_1, \quad\text{ where }\quad\cC_1:= \frac{1}{2\sqrt{2} \Ccont_\star \Creg_\star \Cint n_{\max}},
\eeq
then
\beqs
\N{u-u_h}_{L^2(\OR)} \leq 2
\Ccont_\star \eta(\cH_h) \N{u-w_h}_{\HoDk} \quad \tfa w_h \in \cH_h.
\eeqs
\ele

\bpf
Let $\xi= \cS^*(u-u_h)$; i.e.~$\xi$ is the solution of variational problem
\beqs
\tfind \xi \in \cH \tst a(v,\xi)= (v,u-u_h)_{L^2(\OR)} \quad\tfa v\in \cH. 
\eeqs
Then, by Galerkin orthogonality \eqref{eq:GO2} and the definition of $a_\star(\cdot,\cdot)$ \eqref{eq:astar}, for all $v_h \in \cH_h$,
\begin{align}\label{eq:Schatz_ref}
\N{u-u_h}^2_{L^2(\OR)} &= a( u-u_h, \xi) = a( u-u_h,\xi-v_h ), \\
&= a_\star(u-u_h, \xi-v_h)_{L^2(\OR)} -k^2 (n(u-u_h), \xi-v_h)_{L^2(\OR)}.\nonumber
\end{align}
We choose $v_h= \cP_h\xi$, and then use (in the following order) (i) the Galerkin orthogonality \eqref{eq:GO2}, (ii) continuity of $a_\star(\cdot,\cdot)$, (iii) the bound \eqref{eq:ep2}, (iv) the upper bound in the norm equivalence \eqref{eq:starnormequiv} and the bound \eqref{eq:ep1}, and (v) the consequence \eqref{eq:eta2} of the definition of $\eta$ to obtain that, for all $w_h\in \cH_H$,
\begin{align}\label{eq:epvsSch}
&\N{u-u_h}^2_{L^2(\OR)}  = a_\star(u-w_h, \xi-\cP_h \xi)_{L^2(\OR)} -k^2 (n(u-u_h), \xi- \cP_h\xi)_{L^2(\OR)}\\
\nonumber&\leq \N{u-w_h}_{\star} \N{\xi-\cP_h \xi}_{\star} + k^2 n_{\max}\N{u-u_h}_{L^2(\OR)}\N{\xi-\cP_h \xi}_{L^2(\OR)}\\
\nonumber&\leq \Big( \N{u-w_h}_\star + hk^2 \sqrt{2}\Cint \Creg_\star \sqrt{\Ccont_\star}n_{\max} \N{u-u_h}_{L^2(\OR)}\Big) \N{\xi-\cP_h \xi}_\star\\
\nonumber&
\leq  \Big( \sqrt{\Ccont_\star} \N{u-w_h}_{H^1_k} + hk^2 \sqrt{2} \Cint\Creg_\star \sqrt{\Ccont_\star}n_{\max}\N{u-u_h}_{L^2}\Big)\\
\nonumber & \hspace{6cm}\times\sqrt{\Ccont_\star} \min_{v_h \in \cH_h} \N{\xi-v_h}_{\HoDk}\\
\nonumber
&\leq  \Big( \sqrt{\Ccont_\star} \N{u-w_h}_{H^1_k} + hk^2 \sqrt{2} \Cint\Creg_\star \sqrt{\Ccont_\star}n_{\max}\N{u-u_h}_{L^2}\Big)\\
\nonumber & \hspace{6cm}\times\sqrt{\Ccont_\star}\eta(\cH_h)\N{u-u_h}_{L^2(\OR)};
\end{align}
the result then follows.
\epf

\bre[Advantage of elliptic-projection over standard duality argument]
Comparing \eqref{eq:Schatz_ref} and \eqref{eq:epvsSch} we see the advantage of the elliptic-projection argument over the 
standard duality argument: in \eqref{eq:epvsSch}, Galerkin orthogonality for $a_\star(\cdot,\cdot)$ has allowed us to obtain $u-w_h$ (with $w_h$ arbitrary) as opposed to $u-u_h$ in the first argument of the sesquilinear form on the right-hand side, leading to the bound \eqref{eq:AN2} instead of \eqref{eq:Schatz1}. The price for this is that we have an additional 
$L^2$ inner product on the right-hand side of \eqref{eq:epvsSch}, and controlling this leads to the condition \eqref{eq:C1}.
\ere

Recall that, by the Cauchy--Schwarz inequality and the inequality \eqref{eq:CDtN1}, $a(\cdot,\cdot)$ is continuous, i.e., for all $u,v \in \cH$,
\beq\label{eq:continuity}
\big|a(u,v)\big| \leq \Ccont \N{u}_{\HoDk} \N{v}_{\HoDk},
\eeq
where $\Ccont:= \max\big\{A_{\max}, n_{\max}\big\} + \CTR_1$.

\ble\label{lem:C1}
Assuming that the Galerkin solution $u_h$ to the variational problem \eqref{eq:FEM} exists, if \eqref{eq:C1} holds, 
then 
\beq\label{eq:AN11}
\N{u-u_h}_{\HoDk} \leq \Big(\cC_2 hk + \cC_3 h k^2 \eta(\cH_h)\Big) \N{u}_{\HoDk},
\eeq
where
\beqs
\cC_2:= \frac{\sqrt{2} \Ccont \Cint \Cosc}{A_{\min}}\quad\tand\quad
\cC_3:=  \frac{4\Ccont_\star \Cint \Cosc \sqrt{n_{\max}+ A_{\min}}}{\sqrt{A_{\min}}}.
\eeqs
\ele

\bpf
Since $\DtN$ satisfies the inequality \eqref{eq:CDtN2}, and $\MA$ and $n$ satisfy the inequalities \eqref{eq:Alimits} and \eqref{eq:nlimits}, 
$a(\cdot,\cdot)$ \eqref{eq:sesqui} satisfies the G\aa rding inequality 
\beq\label{eq:Garding}
\Re a(v,v) \geq A_{\min} \N{v}^2_{\HoDk}- k^2 (n_{\max}+ A_{\min}) \N{v}^2_{L^2(\OR)}.
\eeq
Using Galerkin orthogonality \eqref{eq:GO} and continuity of $a(\cdot,\cdot)$ \eqref{eq:continuity}, we find that that \eqref{eq:key1} holds for any $v_h\in \cH_h$.
Using first the inequality \eqref{eq:Cauchy}
with $\alpha=\|u-u_h\|_{\HoDk}$, $\beta=\Ccont \|u-v_h\|_{\HoDk}$, 
$\eps= A_{\min}$, and then Lemma \ref{lem:AN}, we find that if \eqref{eq:C1} holds, then,  for any $v_h\in \cH_h$, 
\begin{align}\nonumber
&\frac{A_{\min}}{2}\N{u-u_h}^2_{\HoDk} \\ \nonumber
& \hspace{0.5cm}\leq \frac{(\Ccont)^2}{2 A_{\min}} \N{u-v_h}^2_{H^1_k(\OR)} + k^2 \big(n_{\max} + A_{\min}\big) \N{u-u_h}^2_{L^2(\OR)}\\
& \hspace{0.5cm}\leq\left[
\frac{(\Ccont)^2}{2 A_{\min}} + 4k^2 \big(n_{\max} + A_{\min}\big) (\Ccont_\star)^2 \big(\eta(\cH_h)\big)^2
\right]\N{u-v_h}^2_{\HoDk},\label{eq:AN4}
\end{align}
By the consequence \eqref{eq:Cint2} of the definition of $\Cint$ and the bound \eqref{eq:Cosci}/\eqref{eq:Cosci2},
\beq
\N{u-I_h u}_{\HoDk} \leq \sqrt{2}h \Cint |u|_{H^2(\OR)}\leq \sqrt{2}hk \Cint \Cosc \N{u}_{\HoDk}.\label{eq:AN5}
\eeq
Choosing $v_h = I_h u$ in \eqref{eq:AN4}, using \eqref{eq:AN5}, taking the square root and using the inequality $\sqrt{a^2+ b^2}\leq a+b$ for all $a,b>0$, we find the result \eqref{eq:AN11}.
\epf

\bpf[Proof of Theorem \ref{thm:main1}]
Under the assumption that the Galerkin solution $u_h$ exists, the fact that the 
bound \eqref{eq:mainbound1} holds under the condition  \eqref{eq:meshthreshold1}  follows from combining Lemma \ref{lem:C1} with the bound \eqref{eq:etabound1} on $\eta$. 
To prove that $u_h$ exists under the condition
\eqref{eq:meshthreshold1}, recall that, since the variational problem
\eqref{eq:FEM} is equivalent to a linear system of equations in a
  finite-dimensional space,
existence of a solution follows from uniqueness.
Suppose that there exists a $\widetilde{u}_h\in \cH_h$ such that 
$a(\widetilde{u}_h,v_h)=0$ for all 
$v_h\in \cH_h$;
to prove uniqueness, we need to show that $\widetilde{u}_h=0$.
Let $\widetilde{u}$ be such that $a(\widetilde{u},v)=0$ for all $v\in \cH$, so that $\widetilde{u}_h$ is the Galerkin approximation to $\widetilde{u}$.
Repeating the argument in the first part of the proof we see that the condition \eqref{eq:meshthreshold1} holds then the 
bound \eqref{eq:mainbound1} holds (with $u$ replaced by $\widetilde{u}$ and $u_h$ replaced by $\widetilde{u}_h$). 
By Lemma \ref{lem:wellposed}, $\widetilde{u}=0$, so \eqref{eq:mainbound1}  implies that $\widetilde{u}_h=0$ and the proof is complete.
\epf

\bpf[Proof of Theorem \ref{thm:main2}]
This is very similar to the proof of Theorem \ref{thm:main1}, except that we use the bound \eqref{eq:etabound2} on $\eta(\cH_h)$ instead of \eqref{eq:etabound1}.
\epf

\section*{Acknowledgements}

We thank Th\'eophile Chaumont-Frelet (INRIA, Nice), Ivan Graham (University of Bath), and particularly Owen Pembery (University of Bath) for useful discussions.
We thank the referees for numerous useful comments and suggestions.
DL and EAS acknowledge support from EPSRC grant EP/1025995/1. JW was partly supported by Simons Foundation grant 631302.

\bibliographystyle{spmpsci}         
\bibliography{biblio_combined_sncwadditions}

\end{document}